%% file: main.tex
\documentclass[journal]{IEEEtran}

\usepackage{amsmath, amssymb, epsfig, graphicx, bm}
\usepackage{cases}
\usepackage{float, subfigure}
\usepackage{color}
\usepackage{amsfonts}
\usepackage{algorithm}
\usepackage{algpseudocode}
\usepackage{cite, url}
\usepackage[all,cmtip]{xy}
\usepackage{color,hyperref}
\definecolor{darkblue}{rgb}{0,0,0.8}
\hypersetup{colorlinks,breaklinks,
	linkcolor=darkblue,urlcolor=darkblue,anchorcolor=darkblue,citecolor=darkblue}

\newtheorem{theorem}{Theorem}
\newtheorem{lemma}{Lemma}
\newtheorem{assumption}{Assumption}

\newtheorem{remark}{Remark}

\input{section/notation}

\title{\LARGE \bf 
	Compressed Gradient Tracking Methods for Decentralized Optimization with Linear Convergence
}

\author{Yiwei Liao, Zhuorui Li,  Kun Huang, and Shi Pu
	\thanks{This work was partially supported  by the Shenzhen
		Research Institute of Big Data (SRIBD) under Fund No. J00120190011 and by the National Natural Science Foundation of China (NSFC) under grant No. 62003287. Yiwei Liao  and Zhuorui Li contributed equally to this work. 
Corresponding author: Shi Pu.}
	\thanks{Yiwei Liao and Zhuorui Li are with Shenzhen Research Institute of Big Data, Shenzhen, China. Kun Huang and Shi Pu are with the School
		of Data Science, Shenzhen Research Institute of Big Data, The Chinese
		University of Hong Kong, Shenzhen, China.
		{\tt\small (emails: lyw@stu.scu.edu.cn, lizhuorui27@gmail.com, kunhuang@link.cuhk.edu.cn, pushi@cuhk.edu.cn)}}%
}

\begin{document}
	\maketitle
	\thispagestyle{empty}
	\pagestyle{empty}

\begin{abstract}
	Communication compression techniques are of growing interests for solving the decentralized optimization problem under limited communication, where the global objective is to minimize the average of local cost functions over a multi-agent network using only local computation and peer-to-peer communication. In this paper, we first propose a novel compressed gradient tracking algorithm (C-GT) that combines gradient tracking technique with communication compression. In particular, C-GT is compatible with a general class of compression operators that unifies both unbiased and biased compressors. We show that C-GT inherits the advantages of gradient tracking-based algorithms and achieves linear convergence rate for strongly convex and smooth objective functions. In the second part of this paper, we propose an error feedback based compressed gradient tracking algorithm (EF-C-GT) to further improve the algorithm efficiency for biased  compression operators. Numerical examples complement the theoretical findings and demonstrate the efficiency and flexibility of the proposed algorithms.
\end{abstract}
\begin{IEEEkeywords}
	Communication compression, decentralized optimization, gradient tracking, linear convergence, error feedback
\end{IEEEkeywords}
		
\section{Introduction}
In this paper, we study the problem of decentralized optimization over a multi-agent network that consists of  $n$ agents. The goal is to collaboratively solve the following optimization problem:
\begin{equation} \label{problem}
	\begin{array}{c}
		\min\limits_{x\in\RR^p}~f(x):=\frac{1}{n}\sum\limits_{i=1}^n f_i(x),
	\end{array}
\end{equation}
where $x$ is the global decision variable, and each agent has a local objective function $f_i: \RR^p\rightarrow \RR$. The agents are connected through a communication network and can only exchange information with their immediate neighbors in the network. Through local computation and local information exchange, they seek a consensual and optimal solution that minimizes the average of all the local cost functions. 
Decentralized optimization is widely applicable when central controllers or servers are not available or preferable, when centralized communication that involves a large amount of data exchange is prohibitively expensive due to limited communication resources, and when privacy preservation is desirable.

Problem \eqref{problem} has  attracted much attention in recent years and has
found a variety of  applications in wireless networks, distributed control of robotic systems, and machine learning, etc  \cite{Cohen2016Distributed,Nedic2018Distributed,Nedic2020Distributed}. 
Early work considered the distributed subgradient descent (DGD) method with a diminishing step-size policy \cite{Nedic2009distributed}.  Under a constant step-size, EXTRA \cite{Shi2015Extra} first achieved linear convergence rate for strongly convex and smooth cost functions by introducing an extra correction term to DGD. Distributed gradient tracking-based methods were later developed in \cite{Xu2015Augmented,Di2016Next,Nedic2017achieving,Qu2018Harnessing,Pu2020distributed}, where the local gradient descent direction in DGD was replaced by an auxiliary variable that is able to track the average gradient of the local objective functions.  As a result, each agent's local iterate is moving in the global descent direction and converges exponentially  to the optimal
solution for strongly
convex and smooth objective functions \cite{Nedic2017achieving,Qu2018Harnessing}.  Compared with EXTRA, gradient tracking-based methods are also suitable for uncoordinated step-sizes \cite{Xu2015Augmented,Nedic2017Geometrically} and possibly asymmetric weight matrices while preserving linear convergence rates.  
Some variants were also proposed to deal with stochastic gradient information and time-varying or directed network topologies. For example, in \cite{Pu2020distributed}, a distributed stochastic gradient tracking method was considered which exhibits comparable performance to a centralized stochastic gradient algorithm.  
Time-varying networks were considered in \cite{Nedic2015Distributed,Nedic2017achieving,Xie2018Distributed,Sun2020Distributed}, and  
more recent development on directed graphs can be found in \cite{Tsianos2012PushSum,Nedic2015Distributed,Xin2018Linear,Pu2021Push,Xin2020General,Pu2020Robust} and the references therein.
In particular, the Push-Pull/AB methods considered by \cite{Xin2018Linear} and \cite{Pu2021Push} used both row stochastic and column stochastic weight matrices to achieve linear convergence rate for strongly convex  and smooth objective functions over general graphs.

In many application scenarios,  it is vital to design  communication-efficient protocols for distributed computation due to limited communication bandwidth and power constraints. Recently, in order to improve system scalability and communication efficiency, researchers have considered a variety of communication compression techniques, such as sparsification and quantization \cite{1-bit-sgd, NIPS2017_qsgd, signSGD-ICML18, Stich2018Sparsified, karimireddy2019Error, mishchenko2019distributed, tang2019doublesqueeze, Liu2020Double,Stich2020Communication,Beznosikov2020Biased,Xu2020Compressed,Liu2020Distributed}. In the centralized setting, these methods were shown to maintain comparable convergence rates \cite{1-bit-sgd, NIPS2017_qsgd, signSGD-ICML18, Stich2018Sparsified, karimireddy2019Error, mishchenko2019distributed, tang2019doublesqueeze, Liu2020Double}. For decentralized optimization, several techniques were  introduced to alleviate compression errors, including difference compression, extrapolation compression \cite{tang_NIPS2018_7992} and compression error compensation \cite{tang2019doublesqueeze, koloskova*2020decentralized}. A novel algorithm with communication compression, which combines with DGD and preserves the model average,  was presented in \cite{ksj2019choco, koloskova*2020decentralized}. But the method converges sublinearly even when the objective functions are strongly convex. In \cite{Liu2020Linear}, a linearly convergent decentralized optimization algorithm  with compression (LEAD) was introduced for strongly convex and smooth objective functions. The method was based on NIDS \cite{li2019decentralized}, a sibling of EXTRA. 

In light of the advantages of gradient tracking-based methods for decentralized optimization, it is natural to consider the marriage between gradient tracking and communication compression. The first such effort was made in \cite{Kajiyama2020Linear} which considered a quantized gradient tracking method based on a special quantizer. It was shown to achieve linear convergence rate for strongly convex and smooth objective functions. However, the algorithm design is rather complicated and relies on a specific quantizer. In addition, the convergence conditions are not easy to verify. 

In this paper, we first consider a novel gradient tracking-based method (C-GT) for decentralized optimization with communication compression. The algorithm compresses both the decision variables and the gradient trackers to provide a communication-efficient implementation. Unlike the existing methods which are mostly based on unbiased compressors or biased but contractive compressors, C-GT is provably efficient for a general class of compressors, including those which are neither unbiased nor biased but contractive, e.g.,  norm-sign compression methods. We show that C-GT achieves linear convergence for strongly convex and smooth objective functions under such a general class of communication compression techniques. 

In the second part of the paper, we propose an error feedback based compressed gradient tracking algorithm (EF-C-GT) to further improve the algorithm efficiency for biased  compression operators in particular.
Compared with unbiased ones, biased compressors show advantages through average capacity of preserving  information \cite{Beznosikov2020Biased} or test accuracy \cite{Vogels2020PowerSGD}. However, the simple distributed gradient descent method may lead to divergence at an exponential rate if the compression operators are allowed to be biased, and a counter-example was provided for Top-1 compressor in \cite{Beznosikov2020Biased}. More discussion and comparison between biased and unbiased compression operators can be found in \cite{Beznosikov2020Biased}. Error feedback is  the  known technique
that can fix the issue and cope with errors induced by biased, contractive compressors \cite{Stich2018Sparsified,karimireddy2019Error,Stich2020Communication,Stich2020ErrorFeedback,Gorbunov2020Linearly,Horvath2020Better}. For example, the type of methods called error compensation or error correction was developed earlier for a particular application in \cite{1-bit-sgd,Strom2015Scalable,Wu2018Error}. The performance of distributed stochastic gradient descent (SGD) with error feedback was analyzed  on the heterogeneous  data in \cite{Stich2020Communication}. 
We show that EF-C-GT also achieves linear convergence rate for strongly convex and smooth objective functions, and has superior performance over C-GT in numerical experiments.  



The main contributions of the paper are summarized as follows:
\begin{itemize}
	\item We propose a novel compressed gradient tracking  algorithm (C-GT)  for decentralized optimization, which inherits the advantages of gradient tracking-based methods and saves communication costs at the same time.
	\item The proposed C-GT algorithm is applicable to a general class of compression operators and  works under arbitrary compression precision. In particular, the general condition on the compression operators unifies the commonly considered unbiased and biased but contractive compressors and also includes other compression methods such as norm-sign compressors.
	\item C-GT provably achieves linear convergence for minimizing strongly convex and smooth objective functions under the general condition on the compression operators.
	\item We propose EF-C-GT to improve the algorithm efficiency for biased compression methods and prove its linearly convergent property under strongly convex and smooth objective functions.
	\item Simulation examples show that C-GT is efficient and widely applicable to various compressors, and EF-C-GT outperforms C-GT for biased compression methods such as Top-k and Random-k.
\end{itemize}

The rest of this paper is organized as follows. 
We present the C-GT algorithm in Section \ref{sec: Alg}. In Section \ref{sec: CA-CGT}, we perform the convergence analysis for C-GT. We introduce EF-C-GT in Section \ref{sec: Alg-EF} and its convergence result in Section \ref{sec:CA-EFCGT}.  Numerical examples are provided in Section \ref{sec: simulation}. Finally, concluding remarks are given in Section \ref{sec: conclusion}.

\subsection{Notation}\label{subsec: Notation}
Vectors are columns if not otherwise specified in this paper.
Let each agent $i$ hold a local copy $\vx_i\in\mathbb{R}^p$ of the decision variable and a gradient tracker (auxiliary variable) $\vy_i\in\mathbb{R}^p$. At the $k$-th iteration, their values are denoted by $\vx_{i}^{k}$ and $\vy_{i}^{k}$, respectively. 
For notational convenience, define
\begin{align*}
	\vX &:= [\vx_1, \vx_2, \ldots, \vx_n]^{\T}\in\mathbb{R}^{n\times p},\\ \vY &:= [\vy_1, \vy_2, \ldots, \vy_n]^{\T}\in\mathbb{R}^{n\times p},
\end{align*}
and
\begin{equation*}
	\oX :=  \frac{1}{n}\mathbf{1}^{\T} \vX\in\mathbb{R}^{1\times p},\ \ \oY :=  \frac{1}{n}\mathbf{1}^{\T}\vY\in\mathbb{R}^{1\times p},
\end{equation*}
where $\mathbf{1}$ is the column vector with each entry given by 1. At the $k$-th iteration, their values are denoted by $\vX^{k}$, $\vY^{k}$, $\oX^{k}$ and $\oY^{k}$, respectively. 
Auxiliary variables of the agents (in an aggregative matrix form) $\vH_{x}$, $\vH_{y}$, $\vQ_{x}$, $\vQ_{y}$, $\widehat{\vX}$, $\widehat{\vY}$, $\widehat{\vQ}_{x}$ and $\widehat{\vQ}_{y}$ are defined similarly.

An aggregate objective function of the local variables is defined as:
\begin{equation}
	F(\vX):=\sum_{i=1}^nf_i(\vx_i).
\end{equation}
Denote 
\begin{equation*}
	\nabla F(\vX):=\left[\nabla f_1(\vx_1), \nabla f_2(\vx_2), \ldots, \nabla f_n(\vx_n)\right]^{\T}\in\mathbb{R}^{n\times p},
\end{equation*}
and 
\begin{equation}\label{eq:avgF}
	{\nabla} \overline {\vF}(\vX) := \frac{1}{n}{\vone^\T \nabla \vF(\vX)}.
\end{equation}

The inner product of  vectors $\va,\vb\in\mathbb{R}^{p}$  is written as $\langle \va,\vb\rangle$. For  matrices $\vA,\vB\in\mathbb{R}^{n\times p}$, we let $\langle \vA,\vB\rangle$ be the Frobenius inner product.
We use $\|\cdot\|$ to denote the Frobenius norm of vectors and matrices by default. Specially, for square matrices, $\|\cdot\|$ represents the spectral norm. The spectral radius of a square matrix $\vM$ is denoted by $\rho(\vM)$.
		
\section{Problem formulation}\label{sec: Pre}
In this section, we provide the assumptions on the communication graphs and the objective functions.  Then, we discuss different kinds of  compression methods and provide a general description for compression operators.

\subsection{Preliminaries}
We start with introducing the conditions on the communication network/graph and the objective functions.
Assume the agents are connected over a directed network  $\mathcal{G}=(\mathcal{V},\mathcal{E})$, where $\mathcal{V}=\{1,2,\ldots,n\}$ is the set of vertices (nodes) and $\mathcal{E}\subseteq \mathcal{V}\times \mathcal{V}$  consists of ordered pairs of vertices. The ordered pair $(i,j)\in \mathcal{E}$ indicates that there is a directed
edge from agent $i$ to agent $j$ and thus the $i$-th agent can directly send information to the $j$-th agent.
For an arbitrary agent $i\in\mathcal{V}$, we define the set of its in-neighbors as $\mathcal{N}_{i}^{\text{in}}=\left\{j\big|(j,i)\in\mathcal{E}\right\}$  and the set of out-neighbors as
$\mathcal{N}_{i}^{\text{out}}=\left\{j\big|(i,j)\in\mathcal{E}\right\}$.   The cardinality of $\mathcal{N}_{i}^{\text{in}}$ and $\mathcal{N}_{i}^{\text{out}}$, denoted by $\text{Deg}_{i}^{\text{in}}$ and $\text{Deg}_{i}^{\text{out}}$, is referred to as agent $i$'s in-neighbor and out-neighbor degree, respectively.
Regarding the network structure, we make the following standing assumption:
\begin{assumption}\label{Assumption: network}
	The directed graph $\mathcal{G}$ is  strongly connected and permits a nonnegative doubly stochastic weight matrix $\vW=[w_{ij}]\in\RR^{n\times n}$. That is, agent  $i$ can receive  information from agent $j$ if and only if $w_{ij}>0$, and $\vW \vone=\vone$  and $\vone ^{\T}\vW=\vone ^{\T}$.
\end{assumption}
\begin{remark}
	It is possible to construct a doubly stochastic weight matrix for any strongly connected directed graph in theory \cite{Gharesifard2012Distributed}. In practice, it is easier to construct such a matrix for certain types of directed graphs, while in general an iterative computation process is needed. For instance, if  $\text{Deg}_{i}^{\text{in}}=\text{Deg}_{i}^{\text{out}}$ for all $i$ (e.g., when $\mathcal{G}$ is undirected), then a doubly stochastic weight matrix can be easily constructed as $\vW=\vI-a\vL$, where $\vI$ is an identity matrix, $\vL$ is the graph Laplacian and $a>0$ is a tuning parameter.
\end{remark}

We make the following assumption on the objective functions.
\begin{assumption}
	\label{Assumption: function}
	The objective function $f$ is $\mu$-strongly convex, and each local cost function $f_i$'s gradient is $L_i$-Lipschitz continuous, i.e., for any $\vx,\vx'\in\mathbb{R}^p$,
	\begin{align}
		& \langle \nabla f(\vx)-\nabla f(\vx'),\vx-\vx'\rangle\ge \mu\|\vx-\vx'\|^2,\\
		& \|\nabla f_i(\vx)-\nabla f_i(\vx')\|\le L_i \|\vx-\vx'\|.
	\end{align}
\end{assumption}
From Assumption \ref{Assumption: function}, the gradient of $f$ is $L$-Lipschitz continuous, where $L=\max{\{L_i\}}$. Note that the problem \eqref{problem} has a unique solution denoted by $\vx^*\in\mathbb{R}^{1\times p}$  under the assumption.

\subsection{Compression Methods}
In this subsection, we introduce some common assumptions on the compression operators and then present a more general and unified assumption.
\subsubsection{Unbiased compression operators}

Denote \textbf{Compress} the compression function. We first consider a general class of unbiased compression methods, in which the variance of the compression error has an upper bound that is linearly proportional to the norm of the variable of interest \cite{NIPS2017_qsgd,wen2017terngrad,mishchenko2019distributed,Liu2020Double,Liu2020Linear}.
\begin{assumption}\label{Assumption:Unbiased}
	The  compression operator $\cQ:\RR^d\rightarrow\RR^d$ associated with function \textbf{Compress} satisfies $\EE \cQ(\vx)=\vx$, and for all $\vx\in\RR^d$, there exists a constant $C\geq0$ such that $\EE\|\cQ(\vx)-\vx \|^2\leq C\|\vx\|^2$.
\end{assumption}
\begin{remark}
	The expectation  is taken with respect to the random vector corresponding to the internal compression randomness of $\cQ$. 
	Some instances of feasible stochastic
	compression operators satisfying Assumption \ref{Assumption:Unbiased} can be found in \cite{Liu2020Linear,Liu2020Double} and the references therein. 
	For example,  we may consider the  unbiased $b$-bits $q$-norm quantization compression method defined as follows: 
	\begin{align}\label{Quant}
		\cQ_{q}(\vx)=\frac{\|\vx\|_{q}} {2^{b-1}} \text{sign}(\vx)  \odot \left \lfloor{\frac{2^{b-1}|\vx|}{\|\vx\|_{q}}+\vu}\right \rfloor,
	\end{align}
	where $\text{sign}()$ is the sign function, $\odot$ is the Hadamard product, $|\vx|$ is the element-wise absolute value of $\vx$, and $\vu$ is a random perturbation vector uniformly distributed in $[0, 1]^{p}$. It has been shown that the compression  operator is unbiased and its compression variance has an upper bound that  is linearly proportional to the norm of the variables \cite{Liu2020Linear}. Note that the agents only need to transmit the norm $\|\vx\|_{q}$, sign$(\vx)$, and integers in the bracket for communication.
\end{remark}

\subsubsection{Biased compression operators}

We also consider the following class of biased compression methods that are common in practice \cite{koloskova*2020decentralized,Stich2020Communication,Beznosikov2020Biased}.
\begin{assumption}\label{Assumption:Biased}
	The  compression operator $\cC_\delta \colon \RR^d \to \RR^d$ associated with function \textbf{Compress} satisfies  
	\begin{align}\label{def:compressor}
		\EE_{\cC_\delta} \norm{\cC_\delta(\vx) -\vx}^2 &\leq (1-\delta)\norm{\vx}^2,& ~\forall \vx \in \RR^d , 
	\end{align}
	where $\delta\in(0,1]$.
\end{assumption}
\begin{remark}
	If $\delta=1$, there is no compression error, i.e., $\cC_\delta(\vx) =\vx$. 
	Below we give two examples of biased compression operators, i.e., Top-k and Random-k, where $\delta$ satisfies $\delta=\frac{k}{p}$ (see e.g., \cite{Beznosikov2020Biased}).
	\begin{itemize}
		\item \textbf{Top-k}: A subset of $\vx$ which corresponds to the $k$ largest absolute values of $\vx$ is chosen, i.e., $\cC_{top}(\vx)=\vx\odot \ve$, where the element of $\ve$ is set to $1$ if the corresponding index is selected, otherwise it is set to $0$. Specifically, we reorder the elements of $\vx$ as
		$|\vx_{i_1}|\geq|\vx_{i_2}|\geq\cdots\geq|\vx_{i_p}|$ and thus $\ve$ satisfies $\ve_{i_l}=1$ for $l\leq k$ and $\ve_{i_l}=0$ for $l>k$.
		
		\item \textbf{Random-k}: A set of randomly  selected $k$ elements of $\vx$ is transmitted, i.e.,  $\cC_{rnd}(\vx)=\vx\odot \ve$, where the elements of $\ve$ satisfy
		\begin{align}
			\ve_i=&\left\{
			\begin{array}{ll}
				1, & \hbox{with~ probability~ $\frac{k}{p}$;} \\
				0, & \hbox{with~ probability~ $1-\frac{k}{p}$.}
			\end{array}
			\right.
		\end{align}

	\end{itemize}
\end{remark}

\subsubsection{General compression operators}
We now present a general assumption on the compression operators, which includes Assumptions \ref{Assumption:Unbiased} and \ref{Assumption:Biased} as special cases.
\begin{assumption}\label{Assumption:General}
	The  compression operator $\cC \colon \RR^d \to \RR^d$ associated with function \textbf{Compress} satisfies  
	\begin{align}\label{def:Ncompressor}
		\EE_{\cC} \norm{\cC(\vx) -\vx}^2 &\leq C\norm{\vx}^2,& ~\forall \vx \in \RR^d , 
	\end{align}
	and the $r$-scaling of $\cC$  satisfies 
	\begin{align}\label{def:contract}
		\EE_{\cC} \norm{\frac{\cC(\vx)}{r} -\vx}^2 &\leq (1-\delta)\norm{\vx}^2,& ~\forall \vx \in \RR^d , 
	\end{align}
	for some constants $\delta\in(0,1]$ and  $r>0$.
\end{assumption}

\begin{remark}
	On one hand, if $C<1$, Assumption \ref{Assumption:General} degenerates to Assumption \ref{Assumption:Biased} by setting $r=1$ and $\delta=1-C$. On the other hand, if $\cC$ is unbiased, i.e., $\EE \cC(\vx)=\vx$, then  Assumption \ref{Assumption:General} degenerates to Assumption  \ref{Assumption:Unbiased} by setting $r=C+1$ and $\delta=\frac{1}{C+1}$. In short, Assumption \ref{Assumption:General} gives a unified description of   unbiased and biased compression operators and thus Assumptions  \ref{Assumption:Unbiased} and \ref{Assumption:Biased} can be regarded as its special cases.
	
	However, there also exist compression operators where $\cC$ is biased and $C\geq 1$ in Assumption \ref{Assumption:General}, that is, they do not satisfy Assumptions \ref{Assumption:Unbiased} and  \ref{Assumption:Biased}.
	For instance, for  norm-sign compression operators, i.e.,
	\begin{align}\label{Comq}
		\cC(\vx)=\|\vx\|_{q} \text{sign}(\vx),
	\end{align}
	it can be verified that $\EE_{\cC} \norm{\cC(\vx) -\vx}^2 \leq (p-2)\|\vx\|_{q}^2+\norm{\vx}^2$ for some $p\geq 2$. To see why \eqref{def:contract} holds, we denote the $r$-scaling $\cC$ as $\cC^r(\vx)=\frac{\cC(\vx)}{r}$ for notational convenience subsequently.  Taking $r=p$, we have
	\begin{align}\label{ComInfRe}
		\cC^p(\vx)=\frac{\|\vx\|_{q}}{p} \text{sign}(\vx).
	\end{align}
Then $\EE_{\cC^p} \norm{\cC^p(\vx) -\vx}^2 \leq \left(1-\frac{\norm{x}_{q}^2}{p\norm{x}_2^2}\right)\norm{\vx}^2$.
	For $q=1,2,\infty$, we give the concrete values of $C$ and $\delta$ in Table \ref{table:normsign}.
	\begin{table}[H]
		\begin{center}
			\begin{tabular}{cccc}
				\hline
				~         &$q=1$         &$q=2$         &$q=\infty$        \\
				\hline
				$C$       &$(p-1)^2$     &$p-1$         &$p-1$             \\
				$\delta$  &$\frac{1}{p}$ &$\frac{1}{p}$ &$\frac{1}{p^2}$   \\
				\hline
			\end{tabular}
		\end{center}
		\caption{Compression coefficients for norm-sign compression methods}
		\label{table:normsign}
	\end{table}
\end{remark}

%
 
 \begin{remark}
 	 Although some compression operators (e.g., norm-sign) can be rescaled so that the new compression operator $\cC^p$ satisfies the contractive condition in Assumption \ref{Assumption:Biased}, applying the rescaled operator may hurt the performance of the algorithm when compared with directly using the original compression operator $\cC$ (see the numerical example in Section \ref{subsec:norm-sign}). Considering Assumption \ref{Assumption:General} provides us with more flexibility in choosing the most suitable compression method.
 \end{remark}

\section{A Compressed Gradient Tracking Algorithm}\label{sec: Alg}
\input{section/alg/algCGT}

\section{Convergence Analysis for C-GT}\label{sec: CA-CGT}

In this section, we study the convergence properties of the
proposed compressed gradient tracking algorithm for minimizing strongly convex and  smooth cost functions. 
Our analysis relies on constructing a linear system of inequalities that is related to the
optimization error $\Omega_{o}^{k}:=\EE\big[\|\oX^{k}-\vx^*\|^2\big]$, 
consensus error $\Omega_{c}^{k}:=\EE\big[\|\vX^{k}-\vone\oX^{k}\|^2\big]$, 
gradient tracking error $\Omega_{g}^{k}:=\EE\big[\|\vY^{k}-\vone\oY^{k}\|^2\big]$, 
and compression errors $\Omega_{cx}^{k}:=\EE\big[\| \vX^{k}- \vH^{k}_x \|^2\big]$ and  $\Omega_{cy}^{k}:=\EE\big[\|\vY^{k}-\vH^{k}_y\|^2\big]$. 

\subsection{Supporting Lemmas}
In order to derive the main results, we introduce some useful lemmas first.
\begin{lemma}\label{lem1}
	Under Assumption \ref{Assumption: function},  for all $k\geq0$,  there holds
	\begin{equation}
		\|\nabla f((\oX^k)^{\T}) - (\oY^k)^{\T}\|\leq\frac{L}{\sqrt{n}}\|\vX^k-\vone\oX^k\|.
	\end{equation}
	In addition, if $\eta<2/(\mu+L),$ then we have
	\begin{equation}
		\|\vx-\eta\nabla f(\vx) - (\vx^*)^{\T}\|\leq(1-\eta\mu)\|\vx-(\vx^*)^{\T}\|,\forall \vx\in \RR^p.
	\end{equation}
\end{lemma}
\begin{lemma}\label{lem2}
	Suppose Assumption \ref{Assumption: network} holds. Let $\rho_w$ denote the spectral norm of the matrix $\vW-\frac{1}{n}\vone\vone^\T$. Then, for all $\omega\in\RR^{n\times p}$, we have $\rho_w<1$ and 
	\[
	\|\vW\omega-\vone\bar{\omega}\|\leq\rho_w\|\omega-\vone\bar{\omega}\|,
	\]
	where $\bar{\omega}=\frac{1}{n}\vone^\T\omega$.
\end{lemma}

Proofs of Lemma \ref{lem1} and Lemma \ref{lem2} can be found in  Lemma 10 of \cite{Qu2018Harnessing}.
\begin{remark}
	Given any $\gamma \leq 1$, let $\widetilde{\vW}=(1-\gamma)\vI+\gamma \vW$ and denote 
	\begin{equation}
		s := 1-\rho_w,\; \tilde{\rho} := 1-\gamma s.
	\end{equation} 
	We have for all $\omega\in\RR^{n\times p}$ that
	\[
	\|\tW\omega-\vone\bar{\omega}\|\leq\tilde{\rho}\|\omega-\vone\bar{\omega}\|.
	\] 
\end{remark}

The following lemma presents a linear system of inequalities, which is central in our convergence analysis.

\begin{lemma}\label{Lem:CGT}
	Suppose Assumptions \ref{Assumption: network}--\ref{Assumption: function} and \ref{Assumption:General} hold and  $\eta<\min\big\{\frac{2}{\mu+L},\frac{1}{3\mu}\big\}$. We have the following inequalities for Algorithm \ref{Alg:CGT}:
\begin{align}
	\Omega_{o}^{k+1}\leq&(1-\frac{3}{2}\eta\mu)\Omega_{o}^{k} 
	+ \frac{3\eta L^2}{\mu n}\Omega_{c}^{k},\\
	\Omega_{c}^{k+1}\leq&\frac{1+\tilde{\rho}^2}{2}\Omega_{c}^{k} 
	+ c_1\frac{\eta^2}{\gamma} \Omega_{g}^{k}
	+c_2\gamma\Omega_{cx}^{k},\\
\nonumber	\Omega_{g}^{k+1}\leq& nc_3L^2\frac{\eta^2}{\gamma} \Omega_{o}^{k} +\left(c_4L^2\gamma+c_3L^2\frac{\eta^2}{\gamma}\right) \Omega_{c}^{k}\\
	&+ \left(\frac{1 + \tilde{\rho}^2}{2}+\frac{1}{2}c_3\frac{\eta^2}{\gamma}\right)\Omega_{g}^{k}
	+ 3c_2L^2\gamma\Omega_{cx}^{k}
	+ c_2\gamma^2\Omega_{cy}^{k},\\
\nonumber	\Omega_{cx}^{k+1}\leq& 2nt_xL^2\eta^2\Omega_{o}^{k} 
	+ (c_5 \gamma^2+ 2t_xL^2\eta^2)\Omega_{c}^{k}\\
	&+ t_x\eta^2\Omega_{g}^{k}
	+(c_x+c_6\gamma)\Omega_{cx}^{k},\\
\nonumber	\Omega_{cy}^{k+1}\leq& 6nt_yL^4\eta^2\Omega_{o}^{k} 
	+ (4c_7L^{2}\gamma^2 + c_6L^4\eta^2) \Omega_{c}^{k}\\
\nonumber	&+ (3c_8\gamma^2+3t_yL^2\eta^2 )\Omega_{g}^{k} 
	+ 3c_7L^2\gamma^2\Omega_{cx}^{k}\\
	&+ (c_y+c_7\gamma^2)\Omega_{cy}^{k},
\end{align}
	where  positive constants $c_1$--$c_8$ and $c_x,c_y$, $t_x,t_y$ are  given in  Appendix \ref{Pf:LemCGT}.
\end{lemma}

\begin{IEEEproof}
	See Appendix \ref{Pf:LemCGT}
\end{IEEEproof}

\begin{remark}    \label{Rmk}
	When $\eta$ and $\gamma$ are taken as in Lemma \ref{Lem:CGT}, we have the following linear system of inequalities:
	$$
	\vw^{k+1}\leq \vA\vw^k,
	$$
	where 
	$
	\vw^k :=	\big[\Omega_{o}^{k}, 
	\Omega_{c}^{k}, 
	\Omega_{g}^{k}, 
	\Omega_{cx}^{k}, 
	\Omega_{cy}^{k}\big]^{\T}.$ 
The inequality is to be taken component-wisely, and the elements of the transition matrix $\vA = [a_{ij}]$ are corresponding to the parameters in Lemma~\ref{Lem:CGT}.
In light of the inequalities, the errors $\Omega_{o}^{k},
\Omega_{c}^{k},
\Omega_{g}^{k},
\Omega_{cx}^{k},$ and $
\Omega_{cy}^{k}$ all converge to 0 at the linear rate $\cO(\rho({\vA})^k)$ if the spectral radius of $\rho(\vA)$ satisfies $\rho(\vA)<1.$ The following lemma provides a sufficient condition that ensures $\rho(\vA)<1.$  
\end{remark}

\begin{lemma}
	\label{implma}
	(Corollary 8.1.29 in \cite{Horn2012Matrix}) Let $\vM\in\RR^{m\times m}$ be a nonnegative matrix and $\vv\in\RR^m$ be an element-wisely positive vector. If $\vM\vv\le \theta \vv$, then $\rho(\vM)\le \theta.$
\end{lemma}

\subsection{Main Results}
We introduce the main convergence result for C-GT algorithm under Assumption \ref{Assumption:General} in the following theorem, which demonstrates the linearly convergent property of C-GT for the general compression operators.

\begin{theorem}\label{Thm:CGT}
Under Assumption \ref{Assumption:General}, suppose the scaling parameters are given by  $\alpha_x, \alpha_y \in (0, \frac{1}{r}]$ 
and consensus step-size satisfies
\begin{equation*}
	\gamma\leq \min\Big\{1,
\frac{1-c_x}{m_3}\epsilon_4,
\frac{1-c_y}{m_4}\epsilon_5\Big\},
\end{equation*}
	where $s = 1-\rho_w$,  constants $c_x,c_y$, $m_1$-$m_4$ are given in Appendix \ref{Pf:LemCGT}--\ref{Pf:ThmCGT}, and $\epsilon_1$-$\epsilon_5$ satisfy \eqref{condition_epsilon} (see Appendix \ref{Pf:ThmCGT}). Then for a fixed step-size $\eta$ that satisfies
\begin{equation*}
\eta\leq \min\Big\{
\frac{s\epsilon_2}{4\kappa\epsilon_3}\frac{\gamma}{L},
\frac{s\epsilon_3}{12\kappa(2n\epsilon_1+2\epsilon_2+\epsilon_3)}\frac{\gamma}{L}
\Big\},
\end{equation*}
	the spectral radius of $\vA$, i.e., $\rho\big(\vA\big)$, is less or equal to $1-\frac{1}{2}$$\eta\mu$, and hence the optimization error $\Omega_{o}^{k}$ and the consensus error $\Omega_{c}^{k}$ both converge to 0 at the linear rate $\cO((1-\frac{1}{2}\eta\mu)^k)$.
\end{theorem}

\begin{IEEEproof}
	See Appendix \ref{Pf:ThmCGT}.
\end{IEEEproof}
\begin{remark}
	Regardless of the compression constant $C$, we can find parameters $\alpha_{x}$, $\alpha_{y}$ and $\gamma$ to obtain the linear convergence rate for C-GT with a fixed step-size $\eta$. This implies that C-GT can be used to deal with communication compression with arbitrary precision, i.e., for any $C>0$.
\end{remark}
\begin{remark}
	If Assumption \ref{Assumption:General} degenerates into  Assumption \ref{Assumption:Biased}, i.e., $C=1-\delta<1$, then  we have $\alpha_x, \alpha_y \in (0,1]$ for $r=1$. If the compression operators are unbiased, then Assumption \ref{Assumption:General} degenerates into  Assumption \ref{Assumption:Unbiased} and $r$ is given by $r=1$ for $C<1$ and $r=C+1$ for $C\geq 1$. Thus, $\alpha_x, \alpha_y$ satisfy $\alpha_x, \alpha_y \in (0,1]$ for $C<1$ and $\alpha_x, \alpha_y \in (0,\frac{1}{C+1}]$ for $C\geq 1$ under the unbiased assumption. 
\end{remark}

\begin{remark}
	Note that the parameter settings given in Theorem \ref{Thm:CGT} are only sufficient and relatively conservative. For practical consideration, it is often possible to find better parameters to achieve faster convergence.
\end{remark}

\section{An Error Feedback Based Compressed Gradient Tracking Algorithm}\label{sec: Alg-EF}
\input{section/alg/algEFCGT}

\subsection{Convergence Analysis for EF-C-GT}\label{sec:CA-EFCGT}
We now perform the convergence analysis for EF-C-GT under the same cost functions and Assumption \ref{Assumption:Biased}, that is, we assume the compression operator is biased but contractive. It is worth noting that EF-C-GT also works under unbiased compressors satisfying Assumption \ref{Assumption:Unbiased}, and with a slight modification on the updates, it can deal with biased compression operators that are not contractive.

 Similar to the analysis of C-GT, we still consider constructing a linear system of inequality and introduce two extra terms, i.e., error feedbacks for the decision variable and the gradient tracker compression,  
 $\Omega_{ex}^{k}:=\EE\big[\| \vE_x^{k}\|^2\big]$ and  $\Omega_{ey}^{k}:=\EE\big[\|\vE_y^{k}\|^2\big]$ respectively.
 To derive  the main convergence results for EF-C-GT, we first provide a key lemma as follows.
\begin{lemma}\label{Lem:EFCGT}
	Suppose Assumptions \ref{Assumption: network}--\ref{Assumption: function} and \ref{Assumption:Biased}  hold and  $\eta<\min\big\{\frac{2}{\mu+L},\frac{1}{3\mu}\big\}$. We have the following inequalities for Algorithm \ref{Alg:EFCGT}:
\begin{align}
	\Omega_{o}^{k+1}\leq&(1-\frac{3}{2}\eta\mu)\Omega_{o}^{k} 
	+ \frac{3\eta L^2}{\mu n}\Omega_{c}^{k},\\
	\Omega_{c}^{k+1}\leq&\frac{1+\tilde{\rho}^2}{2}\Omega_{c}^{k} 
	+ d_1\frac{\eta^2}{\gamma} \Omega_{g}^{k}
	+d_2\gamma\Omega_{cx}^{k}
	+\frac{6d_2}{\delta}\gamma\Omega_{ex}^{k},
\end{align}
\begin{align}	
\nonumber	\Omega_{g}^{k+1}\leq& 6nd-1L^2\frac{\eta^2}{\gamma} \Omega_{o}^{k} +\left(3d_2L^2\gamma+6d_1L^4\frac{\eta^2}{\gamma}\right) \Omega_{c}^{k}\\
\nonumber	&+ \left(\frac{1 + \tilde{\rho}^2}{2}+3d_1L^2\frac{\eta^2}{\gamma}\right)\Omega_{g}^{k}
	+ 3d_2L^2\gamma\Omega_{cx}^{k}\\
	&+ d_2\gamma^2\Omega_{cy}^{k}
	+\frac{18d_2}{\delta}L^2\gamma\Omega_{ex}^{k}
	+\frac{6d_2}{\delta}\gamma\Omega_{ey}^{k},\\
\nonumber	\Omega_{cx}^{k+1}\leq& 2nt'_xL^2\eta^2\Omega_{o}^{k} 
	+ (d_3 \gamma^2+ 2t'_xL^2\eta^2)\Omega_{c}^{k}\\
	&+ t'_x\eta^2\Omega_{g}^{k}
	+(d_x+d_3\gamma)\Omega_{cx}^{k}
	+\frac{6d_3}{\delta}\gamma^2\Omega_{ex}^{k},\\
\nonumber	\Omega_{cy}^{k+1}\leq& 6nt'_yL^4\eta^2\Omega_{o}^{k} 
	+ (3t'_yL^{2}\gamma^2 + 6t'_yL^4\eta^2) \Omega_{c}^{k}\\
\nonumber	&+ (d_4\gamma^2+3t'_yL^2\eta^2 )\Omega_{g}^{k} 
	+ 3d_4L^2\gamma^2\Omega_{cx}^{k}\\
	&+ (d_y+d_4\gamma^2)\Omega_{cy}^{k}
	+\frac{18d_4}{\delta}L^2\gamma^2\Omega_{ex}^{k}
	+\frac{6d_4}{\delta}\gamma^2\Omega_{ey}^{k},\\
	\Omega_{ex}^{k+1}\leq& \frac{2(1-\delta)}{\delta}\Omega_{cx}^{k}
	+\left( 1-\frac{\delta}{2}\right)\Omega_{ex}^{k},\\
	\Omega_{ey}^{k+1}\leq& \frac{2(1-\delta)}{\delta}\Omega_{cy}^{k}
	+\left( 1-\frac{\delta}{2}\right)\Omega_{ey}^{k},
\end{align}
\end{lemma}
where positive constants $d_x,d_y$,  $d_1$--$d_4$, and   $t'_x,t'_y$ are given in Appendix \ref{Pf:LemEFCGT}--\ref{Pf:ThmEFCGT}.
\begin{IEEEproof}
	See Appendix \ref{Pf:LemEFCGT}.
\end{IEEEproof}

Based on the above lemma, we present the main convergence result for EF-C-GT in the following theorem, which demonstrates the linearly convergent property of EF-C-GT.
\begin{theorem}\label{Thm:EFCGT}
Under Assumption \ref{Assumption:Biased}, suppose the scaling parameters are given by $\alpha_x, \alpha_y \in (0,1]$
	and the consensus step-size satisfies
\begin{align}	
	\gamma\leq &\min\Big\{1, 
	\frac{1-d_x}{m'_3}\epsilon_4,
	\frac{1-d_y}{m'_4}\epsilon_5\Big\},
\end{align}
where positive constants $d_x,d_y$, $\epsilon_1$--$\epsilon_7$, and $m'_3,m'_4$ are given in Appendix \ref{Pf:LemEFCGT}--\ref{Pf:ThmEFCGT}.  
Then for a fixed step-size $\eta$ that satisfies
\begin{align}
	\eta\le &\min\Big\{
	\frac{s\epsilon_3}{6\kappa(2n\epsilon_1+2\epsilon_2+\epsilon_3)}\frac{\gamma}{L},
	\frac{s\epsilon_2}{4\kappa\epsilon_3}\frac{\gamma}{L},
	\frac{\delta}{2\mu}\Big\},
\end{align}
	the spectral radius of $\vB$ (see \eqref{EFLIS} in Appendix \ref{Pf:ThmEFCGT}), i.e., $\rho\big(\vB\big)$, is less or equal to $1-\frac{1}{2}$$\eta\mu$, and hence the optimization error and the consensus error  both converge to 0 at the linear rate $\cO((1-\frac{1}{2}\eta\mu)^k)$.
\end{theorem}
\begin{IEEEproof}
	See Appendix \ref{Pf:ThmEFCGT}.
\end{IEEEproof}
\section{Numerical Examples}\label{sec: simulation}
	In this part, we provide some numerical examples to confirm our theoretical results and compare with a few different algorithms under various network settings.  
Consider the ridge regression problem:
\begin{align}\label{Ridge Regression}
	\min_{x\in \mathbb{R}^{p}}f(x)=\frac{1}{n}\sum_{i=1}^nf_i(x)\left(=\left(u_i^{\T} x-v_i\right)^2+\rho\|x\|^2\right),
\end{align}
where $\rho>0$ is a penalty parameter. The pair $(u_i,v_i)$ is a sample that belongs to the $i$-th agent, where $u_i\in\mathbb{R}^p$ represents the features and $v_i\in\mathbb{R}$ represents the observations or  outputs.  
In the simulations, pairs $(u_i,v_i)$ are pre-generated: input $u_i\in[-1,1]^p$ is uniformly distributed, and the output $v_i$ satisfies $v_i=u_i^{\T} \tilde{x}_i+\varepsilon_i$, where $\varepsilon_i$ are independent Gaussian noises with mean $0$ and variance $25$, and $\tilde{x}_i$ are predefined parameters evenly located in $[0,1]^p$. 
Then, the $i$-th agent can calculate the gradient of its local objective function $f_i(x)$ with $g_i(x,u_i,v_i)=2(u_i^{\T}x -v_i)u_i+2\rho x$. The unique optimal solution of the problem is $x^*=(\sum_{i=1}^n u_i u_i^{\T}+n\rho\mathbf{I})^{-1}\sum_{i=1}^n u_i v_i$. 

In our experimental settings, we consider penalty parameter $\rho=0.01$. The number of nodes is $n=10$, and the dimension of variables is $p=20$. Meanwhile, $\vx_i^0$ is randomly generated in $[0,1]^p$ and other initial values   satisfy $\vH_{x}^0=\vzero$, $\vH_{y}^0=\vzero$, and $\vY^0 = \nabla \vF(\vX^0)$, $\vE_x^0=\vzero,\vE_y^0=\vzero$.   
For the weight matrix $\vW$, the weights are defined as follows:
\begin{equation*}
	w_{ij}=\begin{cases}
		p_i & \text{if }~j\in \mathcal{N}_{i}^{\text{out}},  \\
		1- \sum_{j\in\mathcal{N}_{i}^{\text{out}}}w_{ij} & \text{if }~j=i,\\
		0 & \text{otherwise},
	\end{cases}
\end{equation*}
where $p_i=0.1$ satisfies $1-\text{Deg}_{i}^{\text{out}}*p_i>0$.

\subsection{Unbiased Compression Operators}
In this case, we consider the unbiased $b$-bits $q$-norm quantization method with  $b=2$ and $q=\infty$ in \eqref{Quant}, since the $\infty$-norm provides the smallest upper bound for the compression variance \cite{Liu2020Linear}.  The scaling parameters and the consensus step-size are all set to 1. 

We first compare C-GT with the known linearly convergent algorithm with communication compression, LEAD \cite{Liu2020Linear}, for decentralized optimization over fixed undirected networks. 

Note that undirected graphs are special cases of directed graphs that satisfies $(i,j)\in\mathcal{E}$ if and only if $(j,i)\in\mathcal{E}$. In particular, we consider a ring network topology. 
The step-sizes are set to
$\eta=0.09$ for C-GT and $0.12$ for LEAD, which are optimal respectively.

\begin{figure}[htp]
	\begin{center}
		\vspace{-1em}
		\includegraphics[width=7cm]{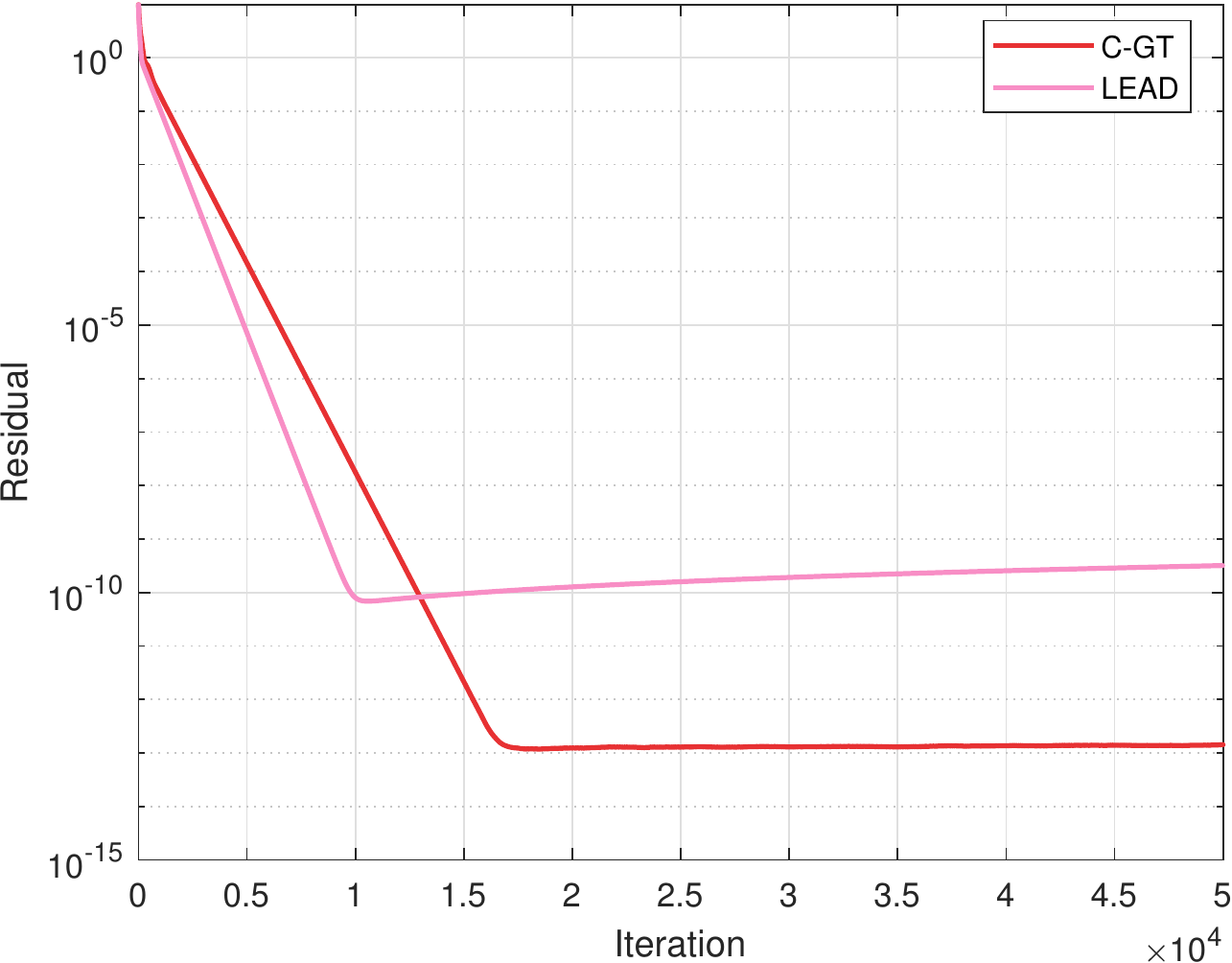}\vspace{-1em}
		\caption{Residuals against number of iterations over an undirected ring graph.}\label{Fig1}
	\end{center}
\end{figure}

It can be seen from Fig. \ref{Fig1} that the optimization errors for both algorithms decrease exponentially fast, which verifies the linearly convergent property of the competing algorithms. Meanwhile, C-GT converges slightly slower than LEAD, but achieves a smaller final error. The two algorithms can be considered comparable in such a case.

\begin{figure}[htp]
	\begin{center}
		\vspace{-1em}
		\includegraphics[width=7cm]{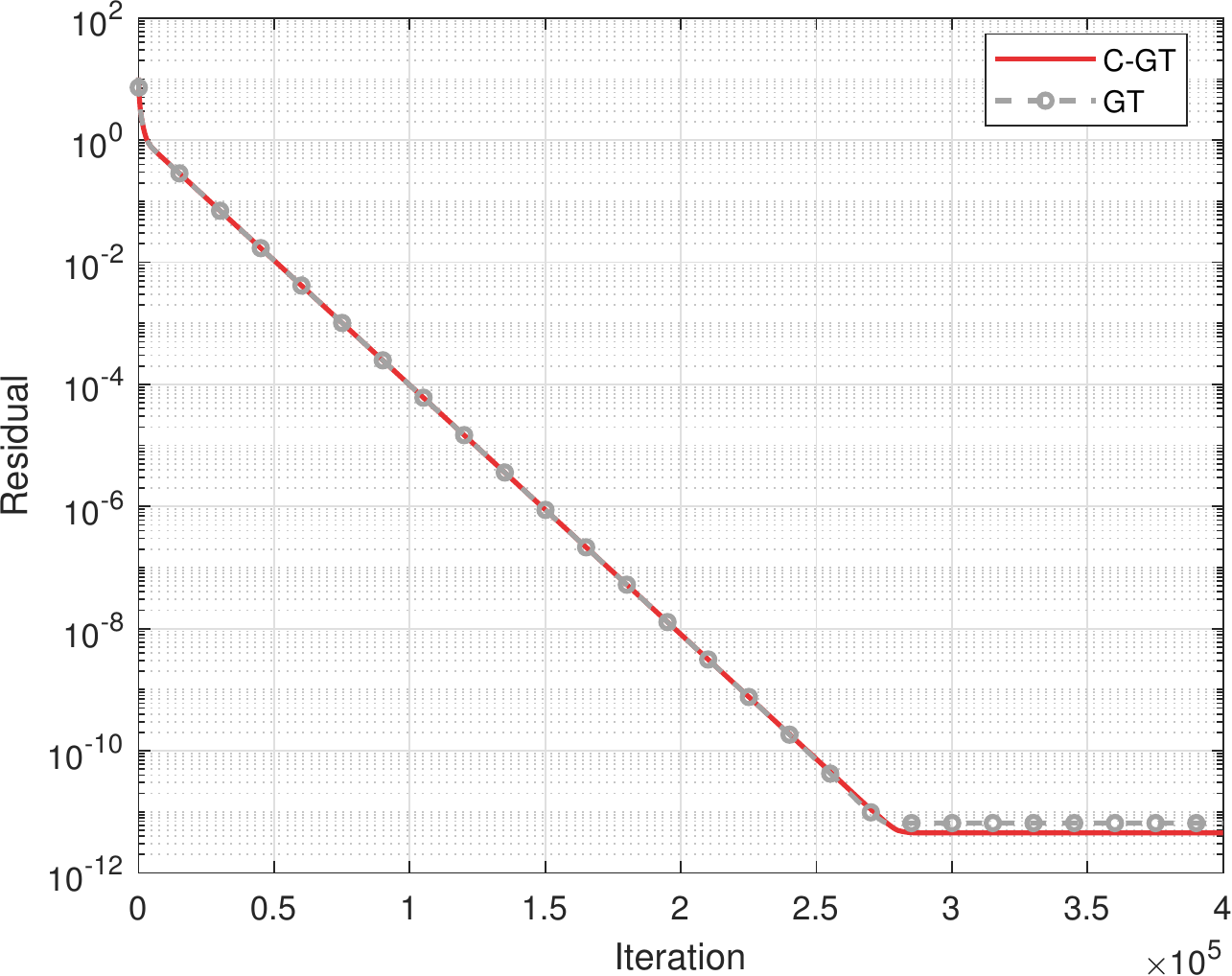}\vspace{-1em}
		\caption{Residuals against number of iterations over  a directed ring graph.}\label{Fig2}
	\end{center}
\end{figure}

Then we consider a directed ring network, for which C-GT is the only known applicable algorithm that achieves linear convergence with communication compression. 
By setting the step-size to $\eta=0.0047$ for the directed ring graph,  
we show in Fig. \ref{Fig2} that C-GT still converges linearly to the optimal solution and achieve almost the same performance as the gradient tracking algorithm without communication compression (GT). These results further demonstrates the effectiveness and flexibility of C-GT.

\subsection{Biased, Contractive Compression Operators}
In this subsection, we consider two sparsification compression methods, i.e., Top-k and Random-k.  We provide simulation results for an undirected network and a directed ring network, respectively. 

\begin{figure}[htp]
	\begin{center}
		\vspace{-1em}
		\subfigure[Undirected ring graph.]{\includegraphics[width=7cm]{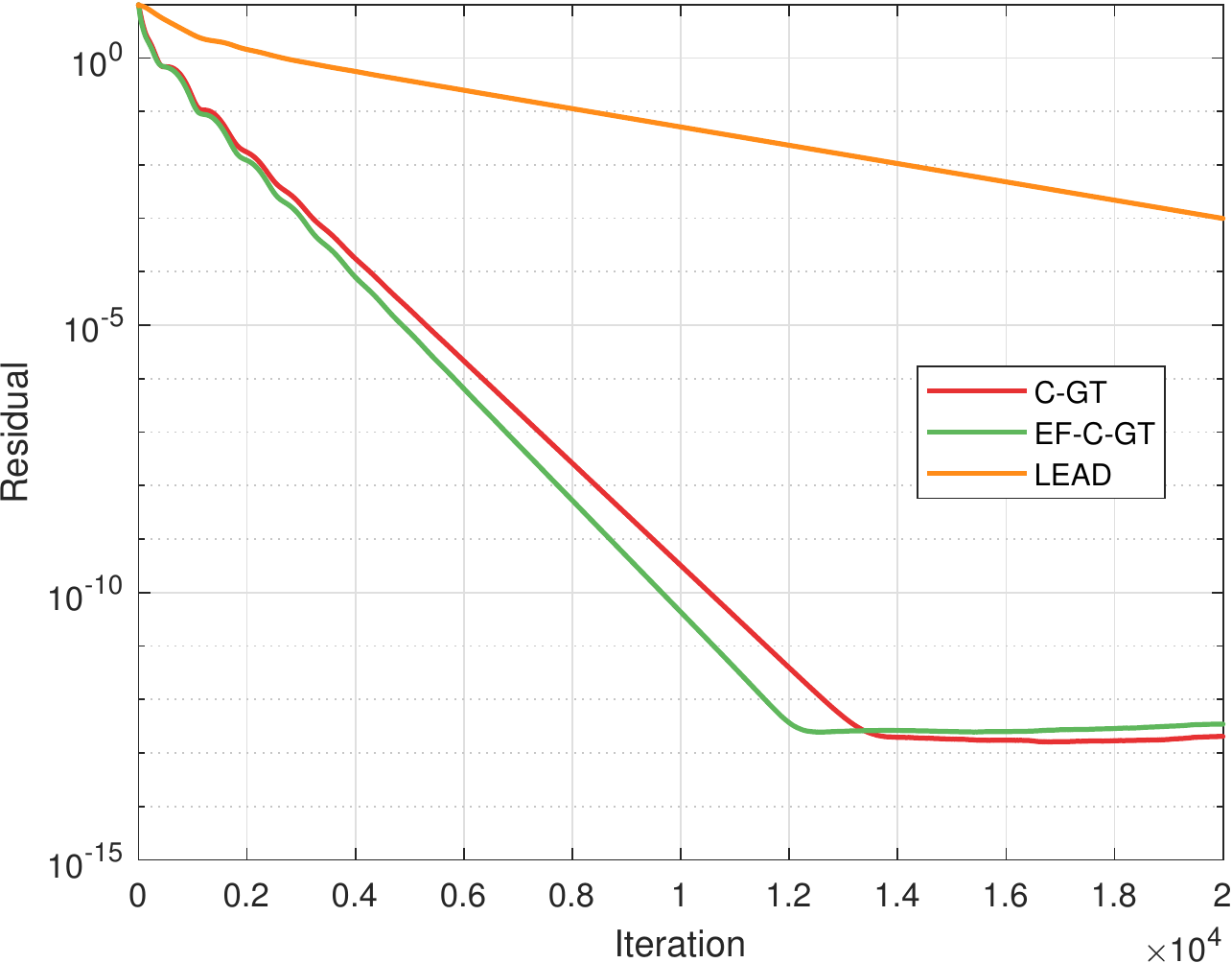}\vspace{-1em}}
		\subfigure[Directed ring graph.]{\includegraphics[width=7cm]{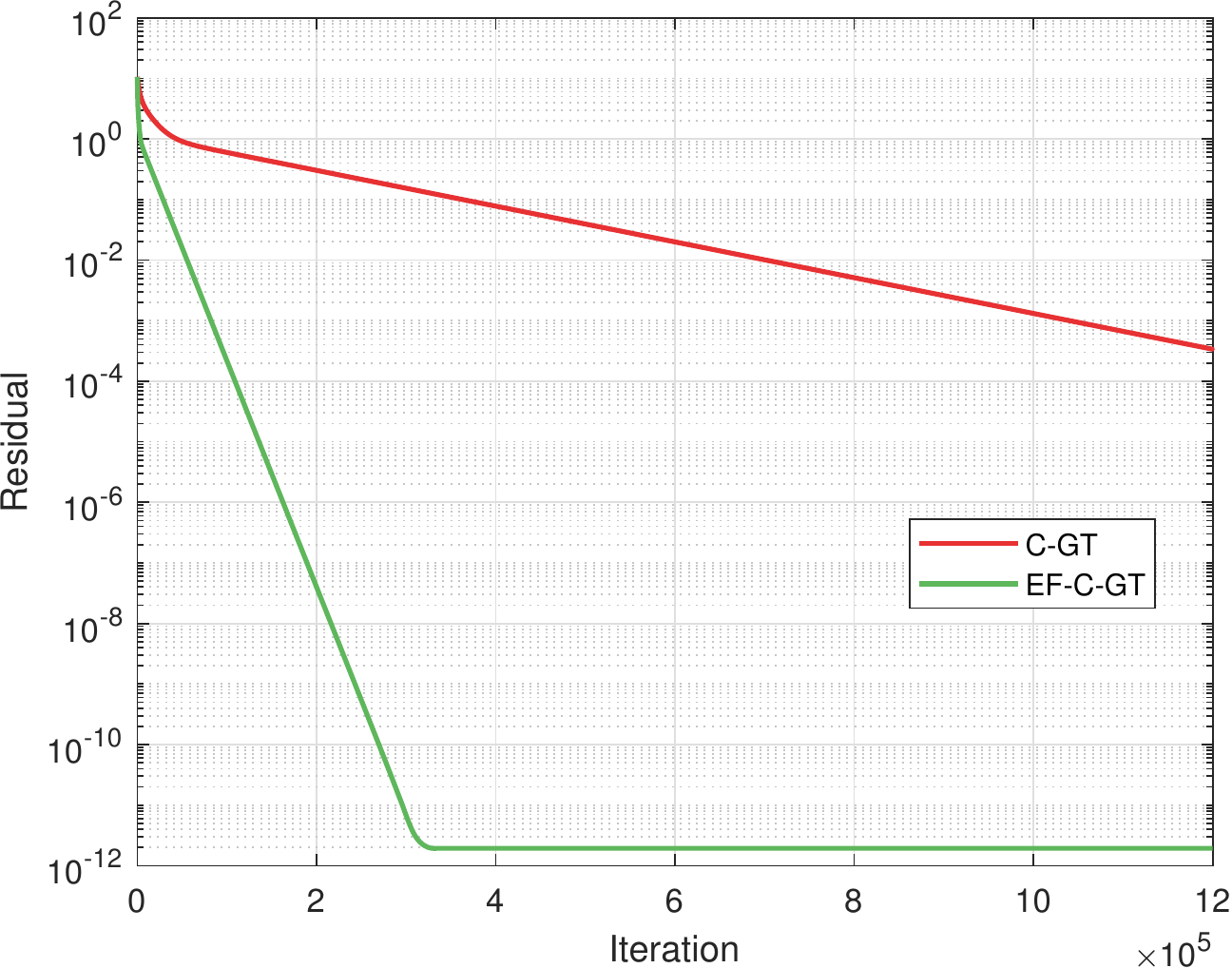}\vspace{-1em}}
		\caption{Residuals against number of iterations for Top-1 compression.}	\label{Fig3}
	\end{center}
\end{figure}

\begin{figure}[htp]	
	\begin{center}
		\vspace{-1em}
		\subfigure[Undirected ring graph.]{\includegraphics[width=7cm]{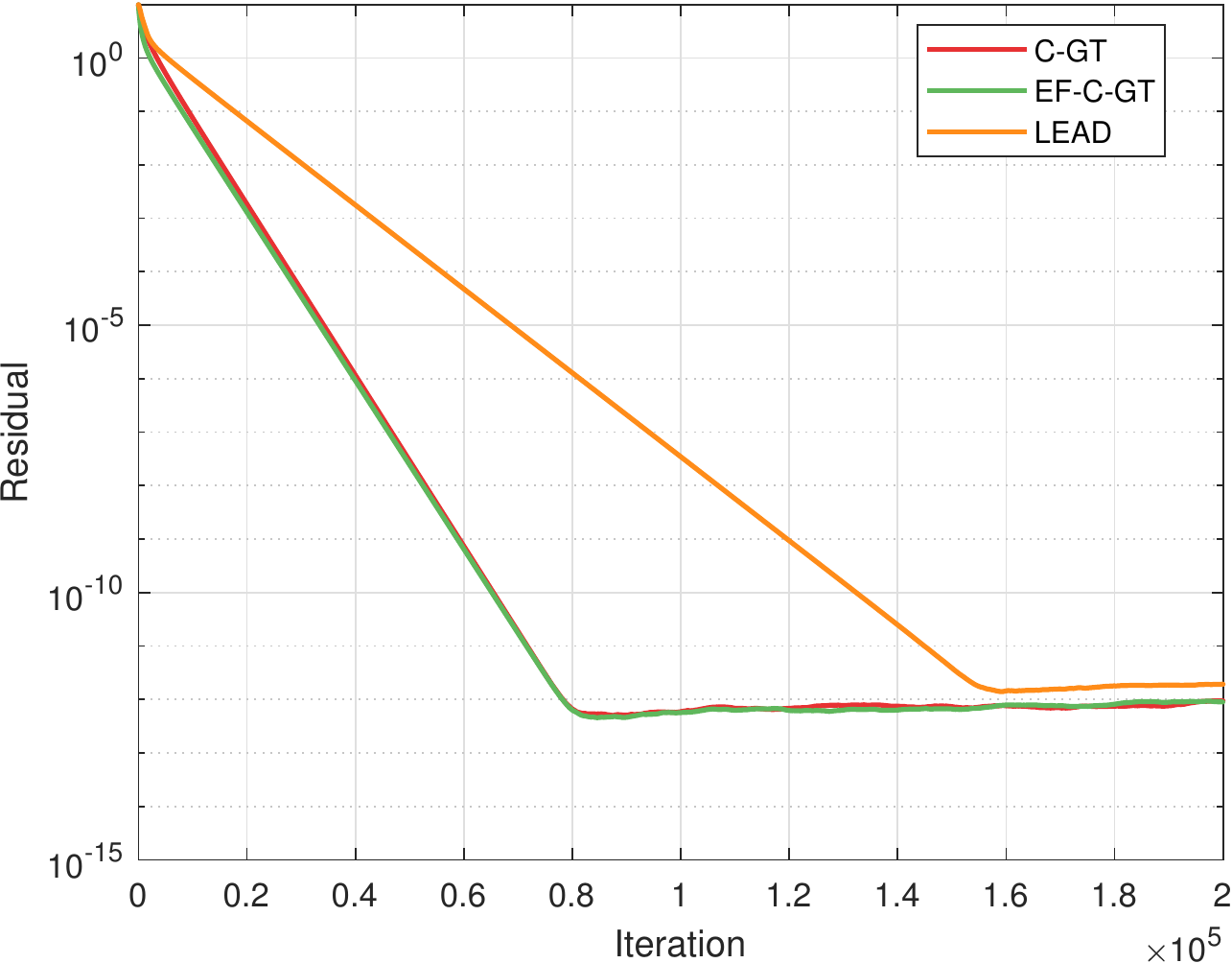}\vspace{-1em}}
		\subfigure[Directed ring graph.]{\includegraphics[width=7cm]{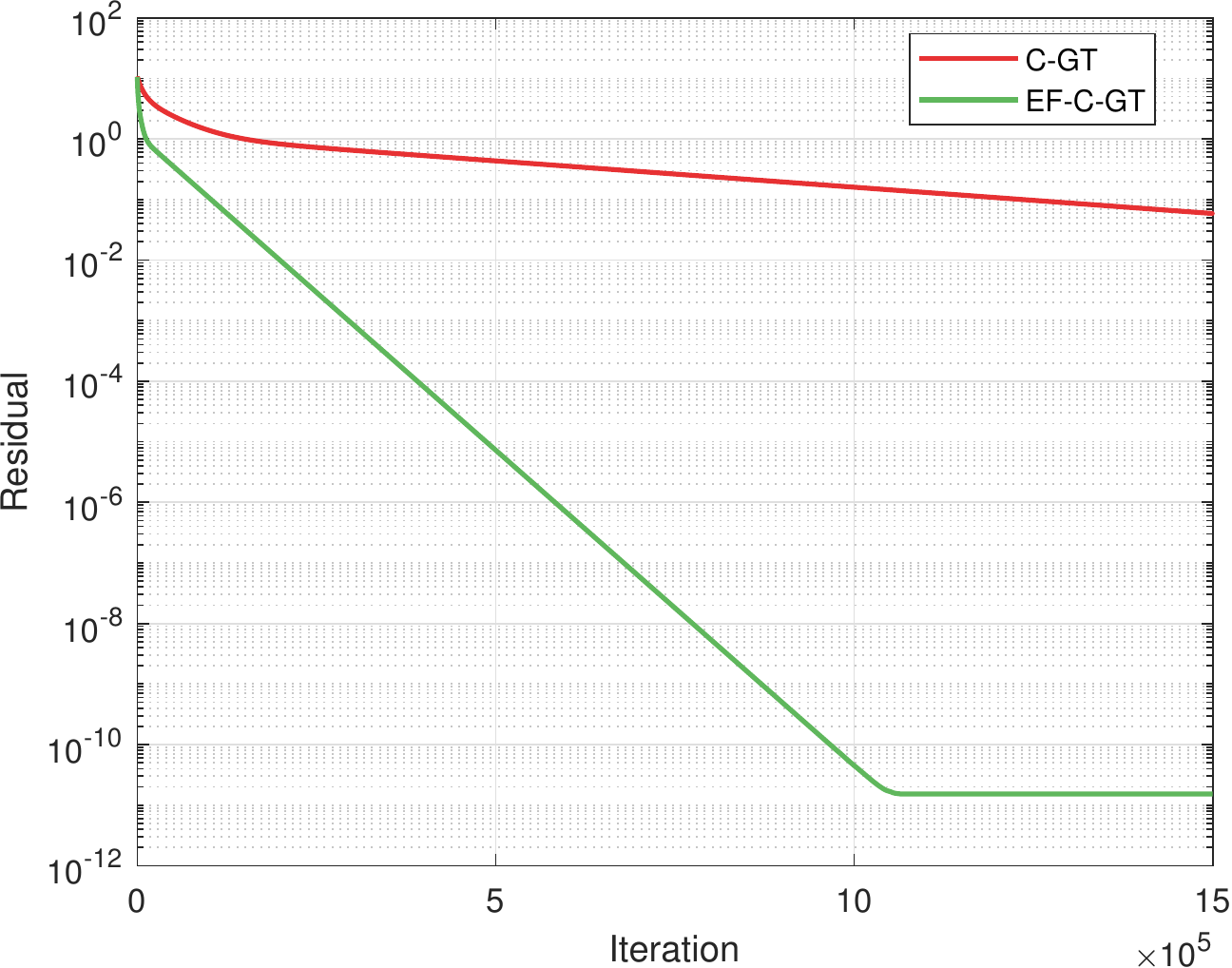}\vspace{-1em}}
		\caption{Residuals against number of iterations for Random-1 compression.}\label{Fig4}
	\end{center}
\end{figure}

For  an undirected ring network, the parameter settings are presented in Table \ref{table:parametersetUnd}. From Fig. \ref{Fig3}(a) and Fig. \ref{Fig4}(a), we can see that C-GT and EF-C-GT both outperform LEAD for both compressors. In addition, EF-C-GT performs slightly better than C-GT for Top-1 compression and similarly to C-GT for Random-1 compression. For a directed ring network where the corresponding parameter settings are given in Table \ref{table:parametersetB},
we can see from Fig. \ref{Fig3}(b)  and Fig. \ref{Fig4}(b)  that EF-C-GT greatly outperforms C-GT for both compressors, which demonstrates the advantages of using error feedback under biased compression operators.

To understand why C-GT performs as well as EF-C-GT for undirected networks, recall the discussion in Section \ref{sec:C-GTalg} where we mention that each agent transmits its total compression error to its neighboring agents and compensate the error locally. Due to the symmetry of undirected networks, such a ``network-based" error compensation is efficient, and thus an additional error feedback term may not be  necessary. On the contrary, for a highly asymmetric directed graph, the error feedback term results in a significant improvement on the algorithm performance.

\begin{table}[htp]
	\begin{center}
		\begin{tabular}{cccccc}
			\hline
			Algorithm       &Compressor   &$\alpha_x$&$\alpha_y$  &$\gamma$       &$\eta$\\
			\hline
			LEAD            &Top-1        &1         &1   &$1\times 10^{-3}$      &0.12\\
			C-GT            &Top-1        &1         &1           &0.6            &0.11\\
			EF-C-GT         &Top-1        &1         &1           &0.6            &0.12\\
			LEAD            &Random-1     &1         &1   &$5\times 10^{-4}$      &0.12\\
		    C-GT            &Random-1     &1         &1           &0.1            &0.11\\
			EF-C-GT         &Random-1     &1         &1           &0.1            &0.11\\	
			\hline
		\end{tabular}
	\end{center}
	\caption{Parameter settings for different algorithms and compression methods  under an undirected ring network.}	
	\label{table:parametersetUnd}
\end{table}
\begin{table}[htp]
	\begin{center}
		\begin{tabular}{cccccc}
			\hline
			Algorithm       &Compressor   &$\alpha_x$&$\alpha_y$  &$\gamma$  &$\eta$     \\
			\hline
			C-GT            &Top-1        &1         &1           &0.5       &0.00034    \\
			EF-C-GT         &Top-1        &1         &1           &1         &0.0043     \\
			C-GT            &Random-1     &1         &1           &0.2       &0.0001     \\
			EF-C-GT         &Random-1     &1         &1           &0.3       &0.0012     \\	
			\hline
		\end{tabular}
	\end{center}
	\caption{Parameter setting  for different algorithms and compression methods under a directed ring network.}	
	\label{table:parametersetB}
\end{table}

It is worth noting that Random-1 compression method amounts to randomized coordinate descent (RCD), and therefore RCD is compatible with the proposed methods.  
Meanwhile, Top-1 compression method results in a greedy coordinate algorithm which is also suitable under our setting.

\subsection{Biased, Non-Contractive Compression Operators}
\label{subsec:norm-sign}
We now consider the  biased norm-sign compression operator \eqref{Comq} and the rescaled norm-sign compressor \eqref{ComInfRe} with $q=\infty$.  
Note that the norm-sign compression operator only satisfies Assumption \ref{Assumption:General}, while its rescaled version satisfies Assumption \ref{Assumption:Biased} by taking $r=p$. 
To make EF-C-GT suitable for the compressors that are not contractive, we make a slight modification on the compression updates: after multiplying   $\vE_x^k$ by  $\beta_x$, we let
$\widehat{\vQ}^{k}_{x}=\textbf{Compress}(\beta_x\vE_{x}^{k}+\vX^{k}-\vH^{k}_{x})$ 
and $\vE_{x}^{k+1}=\beta_x\vE_{x}^{k}+\vX^{k}-\vH^{k}_{x}-\widehat{\vQ}^{k}_{x}$ in Algorithm \ref{Alg:EFCGT}. Similarly, we modify the gradient tracker compression update by multiplying  $\vE_y^k$ with $\beta_y$. The constants $\beta_x$ and $\beta_y$ both belong to $(0,1]$ and should be adjusted according to the constant $C$ in Assumption \ref{Assumption:General}. Indeed, the linear convergence  of the modified algorithm can be  demonstrated similarly to EF-C-GT.

We let $\beta_x=0.01,\beta_y=0.01$, and the other
parameters are provided in Table \ref{table:parametersetG}. 
In the simulation, we use ``N" and ``R" to represent the norm-sign compressor and its rescaled version, respectively. 

\begin{figure}[htp]
	\begin{center}
		\vspace{-1em}
		\includegraphics[width=7cm]{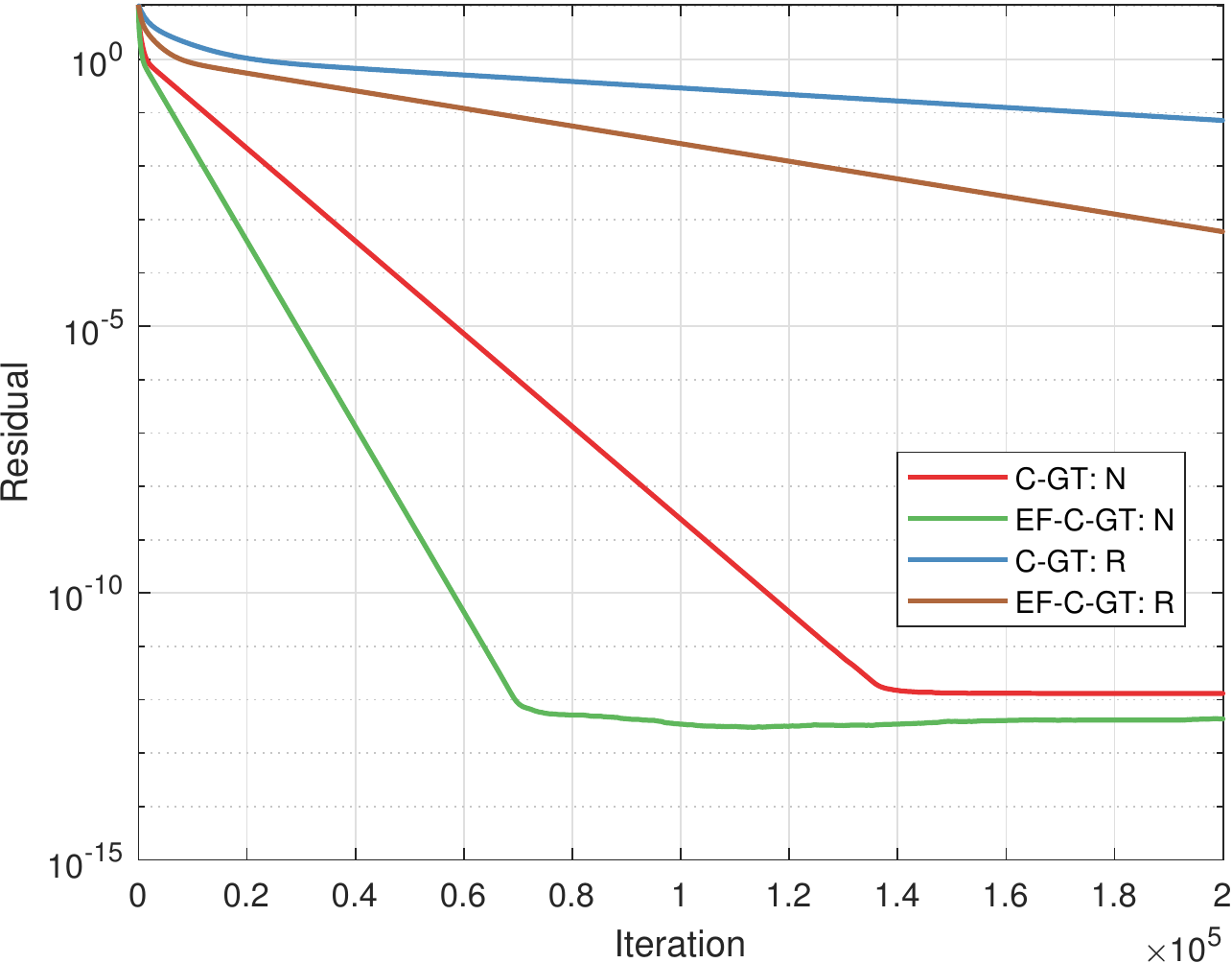}\vspace{-1em}
		\caption{Residuals against number of iterations over  a directed ring graph for the  norm-sign compression operator and its rescaled version.}\label{Fig5}
	\end{center}
\end{figure}

From Fig. \ref{Fig5}, we find that C-GT and EF-C-GT both   work well for the norm-sign compression operator, and EF-C-GT outperforms C-GT. In comparison, using the rescaled compressor leads to slower convergence.  These results suggest that rescaling the compression operator to satisfy the typical contractive requirement (i.e.,  Assumption \ref{Assumption:Biased}) may harm the algorithm performance, and considering Assumption \ref{Assumption:General} provides us with more freedom in choosing the best compression method.

\begin{table}[htp]
	\begin{center}
		\begin{tabular}{cccccc}
			\hline
			Algorithm  &Compressor            &$\alpha_x$&$\alpha_y$  &$\gamma$  &$\eta$\\
			\hline
			C-GT       &Norm-Sign             &0.05      &0.05        &1         &0.01\\
			EF-C-GT    &Norm-Sign             &0.05      &0.05        &1         &0.02\\
			C-GT       &Rescaled Norm-Sign    &1         &1           &0.2       &0.0007\\
			EF-C-GT    &Rescaled Norm-Sign    &1         &1           &0.4       &0.0019\\	
			\hline
		\end{tabular}
	\end{center}
	\caption{Parameter setting for the compression methods}	
	\label{table:parametersetG}
\end{table}

\section{Conclusions}\label{sec: conclusion}
	In this paper, we consider the problem of decentralized optimization with communication compression over a multi-agent network. Specifically, we first propose a compressed gradient tracking algorithm, termed C-GT, and show the algorithm converges linearly for strongly convex and smooth objective functions. C-GT not only inherits the advantages of gradient tracking-based methods, but also works with a wide class of compression operators. To further improve the algorithm efficiency for biased compression methods, we present an error feedback based compressed gradient tracking algorithm (EF-C-GT) and also show its linear convergent property. Simulation examples demonstrate the  effectiveness of C-GT for undirected networks and balanced directed networks, and EF-C-GT outperforms C-GT for some biased compressors such as Top-1 and Random-1.  
	Future work will consider decentralized compression algorithms that work on more general network topologies. We will also consider equipping C-GT with accelerated techniques such as Nesterov's acceleration and momentum methods. Finally, non-convex objective functions are also of future interest.

\section*{ACKNOWLEDGMENT}
	
	We would like to thank Bingyu Wang from The Chinese University of Hong Kong, Shenzhen for providing helpful feedback.
	
\bibliographystyle{IEEEtran}
\bibliography{lywbib/lywbib_compression}	

\appendices

\section{Some Useful Facts}
\input{section/supp/SuppLem}
\section{Proofs for C-GT Algorithm}

\input{section/supp/PfCGT}

\section{Proofs for EF-C-GT Algorithm}
\input{section/supp/PfEFCGT}

\section{Supplementary Algorithms}
\input{section/supp/SuppAlg}

\end{document}

%% file: section/notation.tex
\newcommand{\va}{{\mathbf{a}}}
\newcommand{\vb}{{\mathbf{b}}}

\newcommand{\ve}{{\mathbf{e}}}

\newcommand{\vh}{{\mathbf{h}}}

\newcommand{\vq}{{\mathbf{q}}}

\newcommand{\vu}{{\mathbf{u}}}
\newcommand{\vv}{{\mathbf{v}}}
\newcommand{\vw}{{\mathbf{w}}}
\newcommand{\vx}{{\mathbf{x}}}
\newcommand{\vy}{{\mathbf{y}}}

\newcommand{\vA}{{\mathbf{A}}}
\newcommand{\vB}{{\mathbf{B}}}

\newcommand{\vE}{{\mathbf{E}}}
\newcommand{\vF}{{\mathbf{F}}}

\newcommand{\vH}{{\mathbf{H}}}
\newcommand{\vI}{{\mathbf{I}}}

\newcommand{\vL}{{\mathbf{L}}}
\newcommand{\vM}{{\mathbf{M}}}

\newcommand{\vQ}{{\mathbf{Q}}}

\newcommand{\vU}{{\mathbf{U}}}
\newcommand{\vV}{{\mathbf{V}}}
\newcommand{\vW}{{\mathbf{W}}}
\newcommand{\vX}{{\mathbf{X}}}
\newcommand{\vY}{{\mathbf{Y}}}
\newcommand{\vZ}{{\mathbf{Z}}}

\newcommand{\cC}{{\mathcal{C}}}

\newcommand{\cF}{{\mathcal{F}}}

\newcommand{\cO}{{\mathcal{O}}}

\newcommand{\cQ}{{\mathcal{Q}}}

\newcommand{\EE}{{\mathbb{E}}}
\newcommand{\RR}{\mathbb{R}}
\newcommand{\vzero}{\mathbf{0}}
\newcommand{\vone}{{\mathbf{1}}}
\newcommand{\oX}{\overline{\vX}}
\newcommand{\oY}{\overline{\vY}}
\newcommand{\T}{\intercal}
\newcommand{\tW}{\widetilde{\vW}}
\newcommand{\norm}[1]{\left\| #1 \right\|}      

%% file: section/alg/algCGT.tex
In this section, we introduce the compressed gradient tracking algorithm (C-GT) and its communication-efficient implementation. We also give some interpretations as well as how C-GT connects to existing works.
\subsection{A Compressed Gradient Tracking  Algorithm (C-GT)}\label{sec:C-GTalg}
The proposed compressed gradient tracking  algorithm (C-GT) is presented in Algorithm \ref{Alg:CGT}.\footnote{The implementation details of the algorithms from each agent's view are given in Section \ref{sec:alg_agents_view}.}
\begin{figure}[htp]
	\centering
	\begin{minipage}{.99\linewidth}
		\begin{algorithm}[H]
			\caption{A Compressed Gradient Tracking (C-GT) Algorithm}
			\label{Alg:CGT}
			\textbf{Input:} stopping time $K$, step-size $\eta$, consensus step-size $\gamma$, scaling parameters $\alpha_x, \alpha_y$, and initial values $\vX^0$, $\vH_{x}^0$, $\vH_{y}^0$, $\vY^0 = \nabla \vF(\vX^0)$\\ 
			\noindent\textbf{Output:} $\vX^K, \vY^K$
			\begin{algorithmic}[1]
				\For{$k=0,1,2,\dots, K-1$}
				\State 	$\vQ^{k}_{x}=\textbf{Compress}(\vX^{k}-\vH^{k}_{x})$
				\State $\widehat{\vX}^k = \vH_{x}^k + \vQ_{x}^k$
				\State $\vH^{k+1}_{x}=(1-\alpha_{x})\vH^{k}_{x}+\alpha_{x}\widehat{\vX}^{k}$
				\State $\vQ^{k}_{y}=\textbf{Compress}(\vY^{k}-\vH^{k}_{y})$
				\State $\widehat{\vY}^{k}=\vH^{k}_{y}+\vQ^{k}_{y}$ 
				\State $\vH^{k+1}_{y}=(1-\alpha_{y})\vH^{k}_{y}+\alpha_{y}\widehat{\vY}^{k}$
				\State $\vX^{k+1}=\vX^{k}-\gamma(\vI-\vW)\widehat{\vX}^{k}-\eta \vY^{k}$
				\State $\vY^{k+1}=\vY^{k}-\gamma(\vI-\vW)\widehat{\vY}^{k}+\nabla \vF(\vX^{k+1})-\nabla \vF(\vX^{k})$
				\EndFor
			\end{algorithmic}
		\end{algorithm}
	\end{minipage}
\end{figure}

The function \textbf{Compress} is the compression operator that independently compresses the variables for each agent at every iteration.  In  Line 2,  the difference between $\vX^{k}$ and the auxiliary variable $\vH^{k}_{x}$ is compressed and then added back to $\vH^{k}_{x}$ in Line 3 for obtaining $\widehat{\vX}^k$. Then $\vH^{k+1}_{x}$ is computed as the weighted average of its previous value $\vH^{k}_{x}$ and $\widehat{\vX}^k$ with mixing weight $\alpha_{x}$. Similar updates are conducted for the gradient tracker $\vY^{k}$ in Lines 5-7. 

C-GT performs an implicit error compensation operation that mitigates the impact of 
the compression error, as can be seen from the following argument. The decision variable is updated as
\begin{align}\label{eq:error}
\nonumber 	\vX^{k+1}=&\vX^{k}-\gamma(\vI-\vW)(\vX^{k}-\vE^{k})-\eta \vY^{k}\\
	=&[(1-\gamma)\vI+\gamma \vW] \vX^{k}  -\eta \vY^{k} + \gamma(\vI-\vW) \vE^{k},
\end{align}
where $\vE^{k}:=\vX^{k}-\widehat{\vX}^k$ measures the compression error for the decision variable. The additional  term $(\vI-\vW) \vE^k$ implies that each agent $i$ transmits its total compression error $-\sum_{j\in \mathcal{N}_{i}^{\text{out}}\cup \{i\} } w_{ji}\ve^k_i = -\ve^k_i$ to its neighboring agents and compensates this error locally by adding $\ve^k_i$, where $\ve^k_i\in\mathbb{R}^{1\times p}$ is the  $i$-th row of $\vE^{k}$.  
Similarly, the compression errors for the gradient trackers are also mitigated.

Another key property of C-GT is that gradient tracking is efficient regardless of the compression errors, as for all $k\ge 0$, we have
\begin{align}\label{eq:avgY}
	\nonumber \vone^{\T}\vY^{k+1}=&\vone^{\T}(\vY^{k}-\gamma(\vI-\vW)\widehat{\vY}^{k}+\nabla \vF(\vX^{k+1})-\nabla \vF(\vX^{k}))\\
	\nonumber=&\vone^{\T}\vY^{k}+\vone^{\T}\nabla \vF(\vX^{k+1})-\vone^{\T}\nabla \vF(\vX^{k})\\
	=&\vone^{\T}\nabla \vF(\vX^{k+1}).
\end{align}
The second equality holds because $\vone^{\T}(\vI-\vW)=0$, and the last equality is obtained by induction under the initial condition $\vY^0 = \nabla \vF(\vX^0)$. Therefore, as long as $\vy_{i}^{k}$ reaches (approximate) consensus among all the agents, each $\vy_{i}^{k}$ is able to track the average gradient $\vone^{\T}\nabla \vF(\vX^{k})/n$.
Moreover, by multiplying $\vone^{\T}$ and dividing $n$ on both sides of Line 8, we obtain
\begin{align}
	\label{eq:avgX}\overline \vX^{k+1}=\overline \vX^{k}-\eta \overline{\vY}^{k}
	=\overline \vX^{k}-\eta  {\nabla} \overline \vF(\vX^{k}),
\end{align}
with $ {\nabla} \overline {\vF}(\vX^k)$ defined in \eqref{eq:avgF}. 
Hence the update of $\overline \vX^{k}$ does not contain any compression error.  If all the individual state variables converge to the consensual solution, i.e., $\|\vx_i^k-(\overline\vX^k)^{\T}\| \rightarrow \vzero$,   then we have $\|( {\nabla} \overline{\vF}(\vX^k))^{\T}- {\nabla}f((\overline\vX^k)^{\T})\|\rightarrow \vzero$, and~\eqref{eq:avgX} reduces to the exact gradient descent update.
\begin{remark}
	If no communication compression is performed in the algorithm, i.e., $\widehat{\vX}^k=\vX^{k}$ and $\widehat{\vY}^k=\vY^{k}$,   
	then C-GT recovers the typical gradient tracking algorithm in \cite{Nedic2017achieving,Qu2018Harnessing,Xu2015Augmented}.  To see such a connection, note that C-GT reads
	\begin{align}\label{eq:GT_X}
		\nonumber\vX^{k+1}=&\vX^{k}-\gamma(\vI-\vW)\vX^{k}-\eta \vY^{k}\\
		=&[(1-\gamma)\vI+\gamma \vW] \vX^{k}-\eta \vY^{k},
	\end{align}
	and
	\begin{multline}\label{eq:GT_Y}
\vY^{k+1}=\vY^{k}-\gamma(\vI-\vW)\vY^{k}+\nabla \vF(\vX^{k+1})-\nabla \vF(\vX^{k})\\
		=[(1-\gamma)\vI+\gamma \vW] \vY^{k}+\nabla \vF(\vX^{k+1})-\nabla \vF(\vX^{k}),
	\end{multline}
	where we substitute $\widehat{\vX}^k=\vX^{k}$ and $\widehat{\vY}^k=\vY^{k}$ in Line 8 and Line 9 of C-GT, respectively. By denoting $\widetilde{\vW}:=(1-\gamma)\vI+\gamma \vW$, C-GT takes the same form as the typical gradient tracking method.
\end{remark}

\begin{remark}
	Compared with directly quantizing the decision variables \cite{Kajiyama2020Linear}, the proposed C-GT algorithm reduces the impact of the compression errors through difference compression. It is worth noting that in C-GT, the agents may choose different, uncoordinated step-sizes, which differs from the LEAD algorithm in \cite{Liu2020Linear}.
\end{remark}
\subsection{A Communication Efficient Implementation}
Note that the communication processes associated with computing  $(\vI-\vW) \widehat \vX^k$ and $(\vI-\vW) \widehat \vY^k$ in C-GT require transmitting the full-precision variables $\widehat \vX^k$ and $\widehat \vY^k$ and thus do not enjoy the benefits of compression. In this subsection, we present an equivalent but communication efficient implementation of C-GT for practical consideration, while Algorithm \ref{Alg:CGT} is mainly used for theoretical analysis and  explanation.

\begin{figure}[htp]
	\centering
	\begin{minipage}{.99\linewidth}
		\begin{algorithm}[H]
			\caption{A communication-efficient version of C-GT}
			\label{Alg:CGTe}
			\textbf{Input:} stopping time $K$, step-size $\eta$, consensus step-size $\gamma$, scaling parameters $\alpha_x, \alpha_y$, and initial values $\vX^0$, $\vH_{x}^0$, $\vH_{y}^0$, $\vY^0 = \nabla \vF(\vX^0)$\\ 
			\noindent\textbf{Output:} $\vX^K, \vY^K$
			\begin{algorithmic}[1]
				\State $\vH^0_{x,w} = \vW \vH_{x}^0$
				\State $\vH^0_{y,w} = \vW \vH_{y}^0$				
				\For{$k=0,1,2,\dots, K-1$}
				\State $\widehat{\vX}^k,\widehat{\vX}^k_w,\vH_{x}^{k+1},\vH_{x,w}^{k+1}=\Call{Comm}{\vX^k,\vH^k_x,\vH^k_{x,w}}$\footnote{In  Lines 4 and 5, $\alpha_z$ in the compression function is replaced by $\alpha_x$ for decision difference compression and $\alpha_y$ for gradient difference compression, respectively.}
				\State $\widehat{\vY}^k,\widehat{\vY}^k_w,\vH_{y}^{k+1},\vH_{y,w}^{k+1}=\Call{Comm}{\vY^k,\vH^k_y,\vH^k_{y,w}}$
				\State $\vX^{k+1}=\vX^{k}-\gamma(\widehat{\vX}^{k}-\widehat{\vX}^k_w)-\eta \vY^{k}$
				\State $\vY^{k+1}=\vY^{k}-\gamma(\widehat{\vY}^{k}-\widehat{\vY}^k_w)+\nabla \vF(\vX^{k+1})-\nabla \vF(\vX^{k})$ 
				\EndFor
				\vspace{0.1in}
				\Procedure{Comm}{$\vZ, \vH, \vH_{w}$}
				\State $\vQ =$ \textbf{Compress}($\vZ - \vH$) \hfill $\vartriangleright$ Compression
				\State $\widehat \vZ = \vH + \vQ$
				\State $\widehat \vZ_w = \vH_{w} + \vW\vQ$ \hfill $\vartriangleright$ Communication
				\State $\vH \leftarrow (1-\alpha_z)\vH + \alpha_z \widehat \vZ $
				\State $\vH_{w} \leftarrow (1-\alpha_z)\vH_{w} + \alpha_z \widehat \vZ_w $
				\State \textbf{Return:} $\widehat \vZ, \widehat \vZ_w, \vH, \vH_{w}$
				\EndProcedure
			\end{algorithmic}
		\end{algorithm}
	\end{minipage}
\end{figure}

The main idea of Algorithm~\ref{Alg:CGTe} lies in the procedure $\Call{Comm}{\vZ, \vH, \vH_w}$, in which  two new variables $\vH_w$ and $\widehat\vZ_w$ are introduced. Indeed this procedure is the same as the one considered in \cite{Liu2020Linear}. 
After the  initialization $\vH_{x,w}^0=\vW\vH_x^0$, we obtain $\widehat\vX_w^0=\vH_{x,w}^0+\vW\vQ_x^0=\vW(\vH_x^0+\vQ_x^0)=\vW\widehat\vX^0$ and $\vH_{x,w}^1=(1-\alpha_x)\vH_{x,w}^0+\alpha_x\widehat\vX_w^0
=\vW[(1-\alpha_x)\vH_x^0+\alpha_x \widehat\vX^0]=\vW\vH_x^1$.
Therefore, we ensure $\vH_{x,w}^k=\vW\vH_x^k$ and $\widehat\vX_w^k=\vW\widehat\vX^k$ by induction for all $k$ through the simple communication of the compressed variable $\vQ_x^k$. Similarly, there hold $\vH_{y,w}^k=\vW\vH_y^k$ and $\widehat\vY_w^k=\vW\widehat\vY^k$ for all $k$. 

The procedure $\Call{Comm}{\vZ, \vH, \vH_w}$ works as follows. In Line 10, the difference between $\vZ$ and the auxiliary  variable $\vH$ is compressed into $\vQ$. Then the variable $\vZ$ is estimated based on $\vH$ and the encoded low-bit representation $\vQ$ in Line 11. Through the local communication in Line 12, we obtain $\widehat\vZ_w$, and this is the only required communication step. Finally, given $\widehat\vZ$ and $\widehat\vZ_{w}$, the auxiliary variables $\vH$ and $\vH_{w}$ are updated in Lines 13 and 14, respectively. 

To see why Algorithm~\ref{Alg:CGT} and Algorithm~\ref{Alg:CGTe} are equivalent, note that after the compression and communication procedures  in Lines 4 and 5 respectively, we would have $\widehat\vX_w^k=\vW\widehat\vX^k$ and $\widehat\vY_w^k=\vW\widehat\vY^k$. Thus, the state variable update in Line 6 becomes
\begin{align}\label{eq:EffGT_X}
	\nonumber\vX^{k+1}=&\vX^{k}-\gamma(\widehat{\vX}^{k}-\vW\widehat\vX^k)-\eta \vY^{k}\\
	=&\vX^{k}-\gamma(\vI-\vW)\widehat{\vX}^{k}-\eta \vY^{k},
\end{align}
and the gradient tracker update in Line 7 is given by
\begin{multline}\label{eq:EffGT_Y}
      \vY^{k+1}=\vY^{k}-\gamma(\widehat{\vY}^{k}-\vW\widehat\vY^k)+\nabla \vF(\vX^{k+1})-\nabla \vF(\vX^{k})\\
	=\vY^{k}-\gamma(\vI-\vW)\widehat{\vY}^{k}+\nabla \vF(\vX^{k+1})-\nabla \vF(\vX^{k}).
\end{multline}

\begin{remark}
	In Algorithm~\ref{Alg:CGTe}, computing $\widehat \vX_w^k=\vW \widehat \vX^k$ and $\widehat \vY_w^k=\vW \widehat \vY^k$ does not require the explicit transmission of the full-precision variables $\widehat \vX^k$ and $\widehat \vY^k$, and thus communication efficiency is guaranteed. 
\end{remark}

%% file: section/alg/algEFCGT.tex
In this section, we propose an error feedback based compressed gradient tracking  algorithm (EF-C-GT) to further improve upon the algorithm efficiency of C-GT for biased compression methods particularly.
The algorithm is presented in Algorithm \ref{Alg:EFCGT}. 

Compared to C-GT, in EF-C-GT we consider the difference compression with additional error feedback in Lines 3-4 and 8-9, and $\widehat{\vX}^k,\widehat{\vY}^k$ are updated by $\widehat{\vQ}^k_x$ and $\widehat{\vQ}^k_y$, respectively. The new variables $\vE^k_x,\vE^k_y$ are used to accumulate the compression errors for each node in the network. Then each agent can make use of its history information to correct the bias induced by the biased compressors on the cumulative compression error. The price to pay by considering error feedback lies in twice the number of compression operations in each iteration of the algorithm. 
\begin{figure}[htp]
	\centering
	\begin{minipage}{.99\linewidth}
		\begin{algorithm}[H]
			\caption{An Error Feedback Based C-GT Algorithm (EF-C-GT)}
			\label{Alg:EFCGT}
			\textbf{Input:} stopping time $K$, step-size $\eta$, consensus step-size $\gamma$, scaling parameters $\alpha_x, \alpha_y$, and initial values $\vX^0$, $\vH_{x}^0$, $\vH_{y}^0$, $\vY^0 = \nabla \vF(\vX^0)$\\ 
			\noindent\textbf{Output:} $\vX^K, \vY^K$
			\begin{algorithmic}[1]
				\For{$k=0,1,2,\dots, K-1$}
				\State 	$\vQ^{k}_{x}=\textbf{Compress}(\vX^{k}-\vH^{k}_{x})$
				\State 	$\widehat{\vQ}^{k}_{x}=\textbf{Compress}(\vE_{x}^{k}+\vX^{k}-\vH^{k}_{x})$
				\State 	$\vE_{x}^{k+1}=\vE_{x}^{k}+\vX^{k}-\vH^{k}_{x}-\widehat{\vQ}^{k}_{x}$
				\State $\widehat{\vX}^k = \vH_{x}^k + \widehat{\vQ}^{k}_{x}$
				\State $\vH^{k+1}_{x}=\vH^{k}_{x}+\alpha_{x}\vQ^{k}_{x}$
				\State 	$\vQ^{k}_{y}=\textbf{Compress}(\vY^{k}-\vH^{k}_{y})$
				\State 	$\widehat{\vQ}^{k}_{y}=\textbf{Compress}(\vE_{y}^{k}+\vY^{k}-\vH^{k}_{y})$
				\State 	$\vE_{y}^{k+1}=\vE_{y}^{k}+\vY^{k}-\vH^{k}_{y}-\widehat{\vQ}^{k}_{y}$
				\State $\widehat{\vY}^k = \vH_{y}^k + \widehat{\vQ}^{k}_{y}$
				\State $\vH^{k+1}_{y}=\vH^{k}_{y}+\alpha_{y}\vQ^{k}_{y}$
				\State $\vX^{k+1}=\vX^{k}-\gamma(\vI-\vW)\widehat{\vX}^{k}-\eta \vY^{k}$
				\State $\vY^{k+1}=\vY^{k}-\gamma(\vI-\vW)\widehat{\vY}^{k}+\nabla \vF(\vX^{k+1})-\nabla \vF(\vX^{k})$
				\EndFor
			\end{algorithmic}
		\end{algorithm}
	\end{minipage}
\end{figure}
\begin{remark}
For unbiased compressors, there is no need to consider error feedback since $\EE[\vE_{x}^{k+1}] = 0,\forall k\ge 0$ from Lines 3-4, and the performance of C-GT and EF-C-GT are comparable. In \cite{Horvath2020Better}, it was also shown that error feedback does not provide better performance for the unbiased compressors under the distributed setting.  Indeed, the main idea of considering error feedback in this paper is similar to the usage in \cite{Stich2020Communication}.
\end{remark}

For the sake of completeness, we also present the communication efficient implementation of EF-C-GT in Algorithm \ref{Alg:EFCGTe}.

\begin{figure}[htp]
	\centering
	\begin{minipage}{.99\linewidth}
		\begin{algorithm}[H]
			\caption{A communication-efficient version of EF-C-GT}
			\label{Alg:EFCGTe}
			\textbf{Input:} stopping time $K$, step-size $\eta$, consensus step-size $\gamma$, scaling parameters $\alpha_x, \alpha_y$, and initial values $\vX^0$, $\vH_{x}^0$, $\vH_{y}^0$, $\vY^0 = \nabla \vF(\vX^0)$\\ 
			\noindent\textbf{Output:} $\vX^K, \vY^K$
			\begin{algorithmic}[1]
				\State $\vH^0_{x,w} = \vW \vH_{x}^0$
				\State $\vH^0_{y,w} = \vW \vH_{y}^0$				
				\For{$k=0,1,2,\dots, K-1$}
				\State 
				$\widehat{\vX}^k,\widehat{\vX}^k_w,\vH_{x}^{k+1},\vH_{x,w}^{k+1},\vE_{x}^{k+1}$
				\vfill\hspace{\fill}
				$=\Call{EFComm}{\vX^k,\vH^k_x,\vH^k_{x,w},\vE_{x}^{k}}$
				\State
				$\widehat{\vY}^k,\widehat{\vY}^k_w,\vH_{y}^{k+1},\vH_{y,w}^{k+1},\vE_{y}^{k+1}$
				\vfill\hspace{\fill}
				$=\Call{EFComm}{\vY^k,\vH^k_y,\vH^k_{y,w},\vE_{y}^{k}}$
				\State $\vX^{k+1}=\vX^{k}-\gamma(\widehat{\vX}^{k}-\widehat{\vX}^k_w)-\eta \vY^{k}$
				\State $\vY^{k+1}=\vY^{k}-\gamma(\widehat{\vY}^{k}-\widehat{\vY}^k_w)+\nabla \vF(\vX^{k+1})-\nabla \vF(\vX^{k})$ 
				\EndFor
				\vspace{0.1in}
				\Procedure{EFComm}{$\vZ, \vH, \vH_{w},\vE$}
				\State $\vQ =$ \textbf{Compress}($\vZ - \vH$) \hfill $\vartriangleright$ Compression
				\State $\widehat{\vQ}=$ \textbf{Compress}($\vE+\vZ - \vH$) \hfill $\vartriangleright$ Compression
				\State $\widehat \vZ = \vH + \widehat{\vQ}$
				\State $\widehat \vZ_w = \vH_{w} + \vW\widehat{\vQ}$\hfill $\vartriangleright$ Communication
				\State $\vH \leftarrow \vH + \alpha_z   \vQ $
				\State $\vH_{w} \leftarrow \vH_{w} + \alpha_z  \vW\vQ $\hfill $\vartriangleright$ Communication
				\State $\vE \leftarrow \vE+\vZ-\vH - \widehat \vQ $
				\State \textbf{Return:} $\widehat \vZ, \widehat \vZ_w, \vH, \vH_{w}, \vE$
				\EndProcedure
			\end{algorithmic}
		\end{algorithm}
	\end{minipage}
\end{figure}

%% file: section/supp/SuppLem.tex
\subsection{Vector and Matrix Inequalities}
The following results are often invoked.
\begin{lemma}\label{lem:UV}
	Suppose $\vU,\vV\in \RR^{n\times p},$ then we have the following inequality:
	\begin{align}
	\|\vU+\vV\|^2\leq (1+\tau)\|\vU\|^2 + (1+\frac{1}{\tau})\|\vV\|^2,
    \end{align}
	for any constant $\tau>0$. In particular, taking $\tau=\frac{\delta}{2(1-\delta)}$, we have
$
		(1-\delta)\|\vU+\vV\|^2\leq (1-\frac{\delta}{2})\|\vU\|^2 + \frac{2(1-\delta)}{\delta}\|\vV\|^2.
$
\end{lemma}
\begin{lemma}\label{lem:UV_B}
	For any $\vU,\vV\in \RR^{n\times p}$, the following inequality is satisfied:
	\begin{align}\label{ineq:UV}
		\|\vU+\vV\|^2\leq \tau'\|\vU\|^2 + \frac{\tau'}{\tau'-1}\|\vV\|^2,
	\end{align}
	where $\tau'>1$. 
	In addition, for any $\vU_1,\vU_2,\vU_3 \in \RR^{n\times p}$, we have 
$
		\|\vU_1+\vU_2+\vU_3\|^2\leq \tau'\|\vU_1\|^2 + \frac{2\tau'}{\tau'-1}\left[\|\vU_2\|^2+\|\vU_3\|^2\right]
$ 
	and
$
		\|\vU_1+\vU_2+\vU_3\|^2\leq 3\|\vU_1\|^2 + 3\|\vU_2\|^2+3\|\vU_3\|^2.
$
\end{lemma}

\begin{lemma}
	For any $\vU,\vV\in \RR^{n\times p}$ and $\alpha\in[0,1]$, we have
	\begin{align}\label{ineq:alpha}
		\norm{\alpha \vU+(1-\alpha)\vV}^2\leq \alpha\norm{\vU}^2+(1-\alpha)\norm{\vV}^2.
	\end{align}
\end{lemma}

%% file: section/supp/PfCGT.tex
Let $\cF^k$ be the $\sigma$-algebra generated by $\{\vX^0,\vY^0,\vX^1,\vY^1,\cdots,\vX^{k},\vY^{k}\}$, and define $\EE[ \cdot |\cF^k]$ as the conditional expectation with respect to the compression operator given $\cF^k$.
\subsection{Proof of Lemma \ref{Lem:CGT}}\label{Pf:LemCGT}
Before deriving the linear system of inequality, we bound $\EE\left[\|\vX^k-\widehat{\vX}^k\|^2\middle |\cF^k\right]$ and $\EE\left[\|\vY^k-\widehat{\vY}^k\|^2\middle |\cF^k\right]$, respectively. 
From Lines 3 in Algorithm \ref{Alg:CGT}, we know
\begin{multline}
\EE\left[\norm{\vX^k-\widehat{\vX}^k}^2\middle |\cF^k\right]
= \EE\left[\norm{\vX^k-\vH^k_x-\cC(\vX^k-\vH^k_x)}^2\middle |\cF^k\right]\\
\label{ineqN:Xk_Xkhat}  \leq C\norm{\vX^k-\vH^k_x}^2.
\end{multline}
Similarly, we have
\begin{align}\label{ineqN:Yk_Ykhat}
\EE\left[\norm{\vY^k-\widehat{\vY}^k}^2\middle |\cF^k\right]\leq C\norm{\vY^k-\vH^k_y}^2.
\end{align}
Multiplying $\frac{1}{n}\vone^{\T}$ on both sides of Line 8 in Algorithm \ref{Alg:CGT}, we have
\begin{align}\label{eq:avgX1}
	\overline \vX^{k+1}=\overline \vX^k-\eta \overline{\vY}^k
	=\overline \vX^k-\eta  {\nabla} \overline \vF(\vX^k),
\end{align}

\subsubsection{Optimality error}
For notational  simplicity, denote $\lambda=1-\eta\mu$. According to  \eqref{eq:avgX1},  Lemmas \ref{lem1} and \ref{lem:UV}, we obtain
\begin{align*}
& \norm{\oX^{k+1} - \vx^*}^2 \\
=& \norm{\oX^k-\eta\oY^k-\vx^*}^2\\
=& \norm{(\oX^k-\eta\nabla f((\oX^k)^\T)^\T-\vx^*)+\eta\big(\nabla f((\oX^k)^\T)^\T-\oY^k\big)}^2\\
\leq & (1+\tau_1)\lambda^2 \norm{\oX^k-\vx^*}^2+(1+\frac{1}{\tau_1})\frac{\eta^2 L^2}{n}\norm{\vX^k-\vone \oX^k}^2.
\end{align*}
Taking $\tau_1=\frac{3}{8}\eta\mu$ and noticing that $\eta\leq\frac{1}{3\mu}$, we have
\begin{align*}
&\norm{\oX^{k+1} - \vx^*}^2\\
\leq&(1+\frac{3}{8}\eta\mu)\lambda^2\norm{\oX^k - \vx^*}^2 + (1+\frac{8}{3\eta\mu})\frac{\eta^2 L^2}{n}\norm{\vX^k - \vone \oX^k}^2\\
\leq&(1-\frac{3}{2}\eta\mu) \norm{\oX^k-\vx^*}^2 + \frac{3\eta L^2}{\mu n}\norm{\vX^k-\vone \oX^k}^2.
\end{align*}

\subsubsection{Consensus error}
Recalling Line 8 of Algorithm \ref{Alg:CGT} and relation \eqref{eq:avgX1}, we know
\begin{align*}
&\EE \left[\norm{\vX^{k+1} - \vone \oX^{k+1}}^2\middle|\cF^k\right] \\
=&\EE \left[\norm{\vX^k - \gamma(\vI-\vW)\widehat{\vX}^k -\eta \vY^k-\vone (\oX^k - \eta \oY^k)}^2 \middle| \cF^k\right]\\
=&\EE \bigg[\Big\|\gamma(\vI-\vW)(\vX^k-\widehat{\vX}^k) +\vX^k -\gamma(\vI-\vW)\vX^k -\vone \oX^k\\
&-\eta (\vY^k-\vone \oY^k)\Big\|^2 \bigg| \cF^k\bigg]\\
=&\EE \bigg[\Big\|(\tW\vX^k-\vone\oX^k)-\eta (\vY^k-\vone \oY^k) \\
&+\gamma(\vI-\vW)(\vX^k-\widehat{\vX}^k)\Big\|^2 \bigg| \cF^k \bigg] .
\end{align*}
Based on Lemmas \ref{lem2} and   \ref{lem:UV}, we have 
\begin{align*}
&\EE \left[\norm{\vX^{k+1} - \vone \oX^{k+1}}^2\middle|\cF^k\right] \\
\leq&(1+\tau_2)\tilde{\rho}^2\norm{\vX^k-\vone\oX^k}^2 +2(1+\frac{1}{\tau_2})\bigg(\norm{\eta (\vY^k-\vone \oY^k)}^2\\
&+\EE\bigg[\norm{\gamma(\vI-\vW)(\vX^k-\widehat{\vX}^k)}^2 \bigg| \cF^k\bigg]\bigg).
\end{align*}
Taking $\tau_2 = \frac{1-\tilde{\rho}^2}{2 \tilde{\rho}^2}$, we obtain
\begin{align*}
&\EE \left[\norm{\vX^{k+1} - \vone \oX^{k+1}}^2\middle|\cF^k\right] \\
\leq &\frac{1+\tilde{\rho}^2}{2}\norm{\vX^k-\vone\oX^k}^2 +2\frac{1+\tilde{\rho}^2}{1-\tilde{\rho}^2}\Big(\eta^2\norm{\vY^k-\vone \oY^k}^2\\
&+\gamma^2\norm{\vI-\vW}^2\EE \left[\norm{\vX^k-\widehat{\vX}^k}^2 | \cF^k\right]\Big).
\end{align*}
Let $s=1-\rho_w<1$ and $\gamma s<1$, then $\tilde{\rho}=1-s\gamma$. Meanwhile, we can verify that 
$
2\frac{1+\tilde{\rho}^2}{1-\tilde{\rho}^2}=2\frac{s^2\gamma^2-2\gamma s+2}{\gamma s(2-s\gamma)}<\frac{2}{s\gamma}.
$
Then, we derive
\begin{align}\label{ineq:AvgConsensus}
\nonumber 	&\EE \left[\norm{\vX^{k+1} - \vone \oX^{k+1}}^2\middle|\cF^k\right] \\
\nonumber 	\leq &\frac{1+\tilde{\rho}^2}{2}\norm{\vX^k-\vone\oX^k}^2 +\frac{2\eta^2}{s \gamma}\norm{\vY^k-\vone \oY^k}^2\\
	&+\frac{2\gamma}{s}\norm{\vI-\vW}^2 \EE \left[\norm{\vX^k-\widehat{\vX}^k}^2 \middle| \cF^k\right]
\end{align}
Plugging  \eqref{ineqN:Xk_Xkhat} into \eqref{ineq:AvgConsensus}, we get
\begin{align*}
&\EE \left[\norm{\vX^{k+1} - \vone \oX^{k+1}}^2\middle|\cF^k\right] \\
\leq&\frac{1+\tilde{\rho}^2}{2}\norm{\vX^k-\vone\oX^k}^2 +c_1\frac{\eta^2}{\gamma}\norm{\vY^k-\vone \oY^k}^2\\
&+ c_2\gamma \norm{\vX^k-\vH_x^k}^2,
\end{align*}
where $c_1=\frac{2}{s}$ and $c_2=\frac{2C\gamma}{s}\norm{\vI-\vW}^2$.

\subsubsection{Gradient tracker error}
Firstly, the gradient tracker error is bounded by
\begin{align*}
 &\EE \left[\norm{\vY^{k+1} - \vone \oY^{k+1}}^2\middle|\cF^k\right] \\
=&\EE \left[\norm{(\vY^{k+1} - \vone \oY^k)+(\vone \oY^k-\vone \oY^{k+1})}^2\middle|\cF^k\right] \\
=&\EE \left[\norm{\vY^{k+1} - \vone \oY^k}^2-n\norm{ \oY^k- \oY^{k+1}}^2\middle|\cF^k\right] \\
\leq& \EE \left[\norm{\vY^{k+1}- \vone \oY^k}^2\middle|\cF^k\right] 
\end{align*}
For notational convenience, we use $\nabla^k:=\nabla \vF(\vX^{k})$ instead. Then, from Line 9 in Algorithm \ref{Alg:CGT}, we know
\begin{align*}
&\EE \left[\norm{\vY^{k+1} - \vone \oY^{k+1}}^2\middle|\cF^k\right] \\
\leq& \EE \left[\norm{\vY^{k}-\gamma(\vI-\vW)\widehat{\vY}^k+\nabla^{k+1}-\nabla^k- \vone \oY^k}^2\middle|\cF^k\right] \\
=& \EE \bigg[\Big\|(\tW\vY^k-\vone\oY^k)+\gamma(\vI-\vW)(\vY^{k}-\widehat{\vY}^k)\\
&+(\nabla^{k+1}-\nabla^k)\Big\|^2\bigg|\cF^k\bigg] \\
\leq& (1+\tau_3)\tilde{\rho}^2\norm{\vY^k-\vone\oY^k}^2+2(1+\frac{1}{\tau_3})\EE \bigg[\norm{\nabla^{k+1}-\nabla^k}^2\\
&+\norm{\gamma(\vI-\vW)(\vY^k-\widehat{\vY}^k)}^2 \bigg| \cF^k\bigg]
\end{align*}
Similarly to the derivation of \eqref{ineq:AvgConsensus}, we obtain
\begin{align}
\nonumber &\EE \left[\norm{\vY^{k+1} - \vone \oY^{k+1}}^2\middle|\cF^k\right] \\
\nonumber \leq&\frac{1+\tilde{\rho}^2}{2}\norm{\vY^k-\vone\oY^k}^2 +\frac{2}{s\gamma}\EE\left[\norm{\nabla^{k+1}-\nabla^k}^2\right].\\
\label{YkplusCGT}
&+\frac{2\gamma}{s}\norm{\vI-\vW}^2\EE\left[\norm{\vY^k-\widehat{\vY}^k}^2 \bigg| \cF^k\right]
\end{align}
Next, we derive the bound
\begin{align*}
\EE\left[\norm{\nabla^{k+1}-\nabla^k}^2\middle| \cF^k\right]\leq L^2\EE\left[\norm{\vX^{k+1}-\vX^k}^2\middle| \cF^k\right].
\end{align*}
Thus, we need to bound 
\begin{align}
\nonumber &\EE\left[\norm{\vX^{k+1}-\vX^k}^2\middle| \cF^k\right]\\
\nonumber =&\EE\bigg[\Big\|\gamma(\vI-\vW)(\vX^k-\widehat{\vX}^k)\\
\nonumber &\qquad-\gamma(\vI-\vW)(\vX^k-\vone\oX^k)-\eta\vY^k\Big\|^2\bigg| \cF^k\bigg]\\
\nonumber\leq&3\gamma^2\norm{\vI-\vW}^2\EE\left[\norm{\vX^k-\widehat{\vX}^k}^2\middle| \cF^k\right]\\
\label{ineq:XkplusminusXk}
&+3\gamma^2\norm{\vI-\vW}^2\norm{\vX^k-\vone\oX^k}^2+3\eta^2\norm{\vY^k}^2.
\end{align}
In addition,
\begin{align}
\nonumber \norm{\vY^k}^2=&\norm{\vY^k-\vone\oY^k+\vone\oY^k}^2\\
\label{eq:Yk}=&\norm{\vY^k-\vone\oY^k}^2+n\norm{\oY^k}^2,
\end{align}
and 
\begin{align}
	\nonumber \norm{\oY^k}^2=&\norm{\frac{1}{n}\sum_{i=1}^{n}\nabla f_i((\vx_i^k)^\T)-\nabla f((\oX^k)^\T)+\nabla f((\oX^k)^\T)^\T}^2\\
	\nonumber \leq&2\norm{\frac{1}{n}\sum_{i=1}^{n}[\nabla f_i((\vx_i^k)^\T)-\nabla f((\oX^k)^\T)]}^2\\
	\nonumber &+2\norm{\nabla f((\oX^k)^\T)^\T-\nabla f((\vx^*)^\T)^\T}^2\\
 \leq& \frac{2L^2}{n}\norm{\vX^k-\vone\oX^k}^2+2L^2\norm{\oX^k-\vx^*}^2,
	\label{ineq:Ykbar}
\end{align}
where the last inequality stems from the convexity of 2-norm and Assumption \ref{Assumption: function}.
Plugging \eqref{ineqN:Xk_Xkhat}, \eqref{eq:Yk} and \eqref{ineq:Ykbar} into \eqref{ineq:XkplusminusXk}, we obtain
\begin{multline}
\EE\left[\norm{\vX^{k+1}-\vX^k}^2\middle| \cF^k\right]
\leq 6nL^2\eta^2\norm{\oX^k-\vx^*}^2\\
+(6L^2\eta^2+3\gamma^2\norm{\vI-\vW}^2)\norm{\vX^k-\vone\oX^k}^2
+3\eta^2\norm{\vY^k-\vone\oY^k}^2\\
\label{ineq:XkplusminusXk2}+3C\gamma^2\norm{\vI-\vW}^2\norm{\vX^k-\vH^k_x}^2.
\end{multline}
Therefore, we get
\begin{multline}\label{GradplusminusGrad}
\EE\left[\norm{\nabla^{k+1}-\nabla^k}^2\middle| \cF^k\right]
\leq  6nL^4\eta^2\norm{\oX^k-\vx^*}^2\\
+(6L^4\eta^2+3L^2\gamma^2\norm{\vI-\vW}^2)\norm{\vX^k-\vone\oX^k}^2\\
~~+3L^2\eta^2\norm{\vY^k-\vone\oY^k}^2+3L^2C\gamma^2\norm{\vI-\vW}^2\norm{\vX^k-\vH^k_x}^2.
\end{multline}
Substituting \eqref{GradplusminusGrad} and \eqref{ineqN:Yk_Ykhat} into \eqref{YkplusCGT}, we conclude that 
\begin{align}
\nonumber &\EE \left[\norm{\vY^{k+1} - \vone \oY^{k+1}}^2\middle|\cF^k\right] \\
\nonumber=& c_3nL^2\frac{\eta^2}{\gamma}\norm{\oX^k-\vx^*}^2
\nonumber+(c_3L^2\frac{\eta^2}{\gamma}+c_4L^2\gamma)\norm{\vX^k-\vone\oX^k}^2\\
&\nonumber+(\frac{1+\tilde{\rho}^2}{2}+\frac{1}{2}c_3\frac{\eta^2}{\gamma})\norm{\vY^k-\vone\oY^k}^2
\\
\label{Ykplus1}&+3c_2L^2\gamma\norm{\vX^k-\vH^k_x}^2
+c_2\gamma\norm{\vY^k-\vH_y^k}^2,
\end{align}
where $c_3=\frac{12L^2}{s}$ and $c_4=\frac{6\norm{\vI-\vW}^2}{s}$.

\subsubsection{Difference compression error of decision variables}
From Line 4 in Algorithm \ref{Alg:CGT}, we know
\begin{align}
\nonumber & \norm{\vX^{k+1} - \vH_x^{k+1}}^2\\
\nonumber =& \norm{\vX^{k+1} -\vX^k+\vX^k- \vH_x^k-\alpha_x r \frac{\vQ_x^k}{r}}^2 \\
\nonumber =& \Big\|\vX^{k+1} -\vX^k+\alpha_x r (\vX^k-\vH_x^k-\cC^r(\vX^k-\vH_x^k))\\
\nonumber&+(1-\alpha_x r)(\vX^k- \vH_x^k)\Big\|^2\\
\nonumber  \leq& \tau_x\Big\|\alpha_x r (\vX^k-\vH_x^k-\cC^r(\vX^k-\vH_x^k))\\
\nonumber&+(1-\alpha_x r)(\vX^k- \vH_x^k)\Big\|^2+\frac{\tau_x}{\tau_x-1}\norm{\vX^{k+1} -\vX^k}^2\\
\nonumber \leq& \tau_x\bigg[\alpha_x r\norm{\vX^k-\vH_x^k-\cC^r(\vX^k-\vH_x^k)}^2\\
\label{ineq:HkxError}& +(1-\alpha_x r)\norm{\vX^k- \vH_x^k}^2\bigg]+\frac{\tau_x}{\tau_x-1}\norm{\vX^{k+1} -\vX^k}^2,
\end{align}
where $\alpha_x$ satisfies $0<\alpha_x r\leq 1$.
Taking conditional expectation on both sides of \eqref{ineq:HkxError}, we get
 \begin{align}
 \nonumber & \EE \left[\norm{\vX^{k+1} - \vH_x^{k+1}}^2\middle|\cF^k\right]\\
\nonumber \leq& \tau_x\left[\alpha_x r(1-\delta)+(1-\alpha_x r)\right]\norm{\vX^k-\vH_x^k}^2\\
\label{ineq:EHkxError}&+\frac{\tau_x}{\tau_x-1}\EE \left[\norm{\vX^{k+1} -\vX^k}^2\middle|\cF^k\right],
 \end{align}
where the inequality stems from Assumption \ref{Assumption:General}. 
For simplicity, denote $t_x=\frac{3\tau_x}{\tau_x-1}>1$ and $c_x=\tau_x\left[\alpha_x r(1-\delta)+(1-\alpha_x r)\right]=\tau_x(1-\alpha_x r \delta)<1$. Then, plugging \eqref{ineq:XkplusminusXk2} into \eqref{ineq:EHkxError}, we obtain 
\begin{align}
	\nonumber & \EE \left[\norm{\vX^{k+1} - \vH_x^{k+1}}^2\middle|\cF^k\right]\\
	\nonumber\leq& 2t_xnL^2\eta^2\norm{\oX^k-\vx^*}^2+(2t_xL^2\eta^2+c_5\gamma^2)\norm{\vX^k-\vone\oX^k}^2\\
	&+t_x\eta^2\norm{\vY^k-\vone\oY^k}^2+(c_x+c_6\gamma^2)\norm{\vX^k-\vH^k_x}^2,
\end{align}
where $c_5=t_x \norm{\vI-\vW}^2$ and $c_6=t_xC\norm{\vI-\vW}^2$.
\begin{remark}
	For simplicity, denote $\alpha_r=\alpha_x r$. 
	For $\alpha_r \in (0,1]$ and $\delta\in(0,1)$, we know $0<\alpha_r(1-\delta)+(1-\alpha_r)<1$ and thus we  derive $\tau_x[\alpha_r(1-\delta)+(1-\alpha_r)]<1$ by choosing suitable $1<\tau_x<\frac{1}{\alpha_r(1-\delta)+(1-\alpha_r)}$. For $\alpha_r=1$ and $\delta=1$, we have $\alpha_r(1-\delta)+(1-\alpha_r)=0$ and we can choose an arbitrary  $\tau_x>1$.  Thus, we need $0<\alpha_x\leq \frac{1}{r}$. In fact, the function of $\alpha_r$, i.e., $g_x(\alpha_r)=\alpha_r (1-\delta)+(1-\alpha_r)=1-\delta\alpha_r$, decreases as $\alpha_r$ and thus $g_x(\alpha_r)\geq g_x(1)=1-\delta$. 
\end{remark}
\subsubsection{Difference compression error of gradient trackers}
Based on the same bound method as above, we know
\begin{align}
	\nonumber & \EE \left[\norm{\vY^{k+1} - \vH_y^{k+1}}^2\middle|\cF^k\right]\\
	\nonumber\leq& \tau_y\left[\alpha_y r(1-\delta)+(1-\alpha_y r)\right]\norm{\vY^k-\vH_y^k}^2\\
	\label{ineq:EHkyError}&+\frac{\tau_y}{\tau_y-1}\EE \left[\norm{\vY^{k+1} -\vY^k}^2\middle|\cF^k\right].
\end{align}
In light of \eqref{ineq:EHkyError}, we still need to bound
\begin{align}
\nonumber  &\EE \left[\norm{\vY^{k+1} -\vY^k}^2\middle|\cF^k\right]\\
\nonumber= &\EE \bigg[\Big\|\gamma(\vI-\vW)(\vY^k-\widehat{\vY}^k)-\gamma(\vI-\vW)(\vY^k-\vone\oY^k)\\
\nonumber&+(\nabla^{k+1}-\nabla^k)\Big\|^2\bigg|\cF^k\bigg]\\
\nonumber\leq&
3\gamma^2\norm{\vI-\vW}^2\EE\left[\norm{\vY^k-\widehat{\vY}^k}^2\middle|\cF^k\right]\\
\nonumber &+3\gamma^2\norm{\vI-\vW}^2\norm{\vY^k-\vone\oY^k}^2\\
\label{ineq:EHkyError1}&+3\EE\left[\norm{\nabla^{k+1}-\nabla^k}^2\middle|\cF^k\right].
\end{align}
Plugging \eqref{ineqN:Yk_Ykhat} and \eqref{GradplusminusGrad} into \eqref{ineq:EHkyError1}, we conclude that 
\begin{align}
\nonumber  &\EE \left[\norm{\vY^{k+1} -\vY^k}^2\middle|\cF^k\right]\\
\nonumber \leq&
3C\gamma^2\norm{\vI-\vW}^2\norm{\vY^k-\vH^k_y}^2
+18nL^4\eta^2\norm{\oX^k-\vx^*}^2\\
\nonumber &
+(18L^4\eta^2+9L^2\gamma^2\norm{\vI-\vW}^2)\norm{\vX^k-\vone\oX^k}^2\\
\label{ineq:EHkyError2} \nonumber&+(3\gamma^2\norm{\vI-\vW}^2+9L^2\eta^2)\norm{\vY^k-\vone\oY^k}^2\\
&+9L^2C\gamma^2\norm{\vI-\vW}^2\norm{\vX^k-\vH^k_x}^2.
\end{align}
Substituting \eqref{ineq:EHkyError2} into \eqref{ineq:EHkyError}, we obtain
\begin{align}
\nonumber & \EE \left[\norm{\vY^{k+1} - \vH_y^{k+1}}^2\middle|\cF^k\right]\\
\nonumber  \leq& (c_y+c_7\gamma^2)\norm{\vY^k-\vH_y^k}^2
+6t_ynL^4\eta^2\norm{\oX^k-\vx^*}^2\\
\nonumber &+(6t_yL^4\eta^2+3c_8L^2\gamma^2)\norm{\vX^k-\vone\oX^k}^2\\
\label{ineq:EHkyError4} \nonumber&+(c_8\gamma^2+3t_yL^2\eta^2)\norm{\vY^k-\vone\oY^k}^2\\
&+3c_7L^2 \gamma^2\norm{\vX^k-\vH^k_x}^2,
\end{align}
where $c_y=\tau_y\left[\alpha_y r(1-\delta)+(1-\alpha_y r )\right]=\tau_y(1-\alpha_y r \delta)<1$,
$t_y=\frac{3\tau_y}{\tau_y-1}$, $c_7=t_yC\norm{\vI-\vW}^2$, and $c_8=t_y\norm{\vI-\vW}^2$.


\subsection{Proof of Theorem \ref{Thm:CGT}}\label{Pf:ThmCGT}
In terms of Lemma \ref{Lem:CGT} , we consider the following linear system of inequalities:
\begin{align}\label{LIS}
	\vA\boldsymbol{\epsilon}\leq (1-\frac{1}{2}\eta\mu)\boldsymbol{\epsilon},
\end{align}
where  $\boldsymbol{\epsilon}:=[\epsilon_1,\epsilon_2,L^2\epsilon_3,\epsilon_4,L^2\epsilon_5]^{\T}$, the first three columns of $\vA$ is given by 
\begin{align*}
	\begin{bmatrix}
		1-\frac{3}{2}\eta\mu               &\frac{3\eta L^2}{\mu n}      &0 \\
		0    &\frac{1+\tilde{\rho}^2}{2}   &c_1\frac{\eta^2}{\gamma}     \\
		nc_3L^2\frac{\eta^2}{\gamma}       &c_3L^2\frac{\eta^2}{\gamma}+c_4L^2\gamma
		&\frac{1+\tilde{\rho}^2}{2}+\frac{1}{2}c_3\frac{\eta^2}{\gamma}  \\
		2nt_xL^2\eta^2  &c_5\gamma^2+2t_xL^2\eta^2  &t_x\eta^2  \\
		6nt_yL^4\eta^2  &3c_8L^2\gamma^2+6t_yL^4\eta^2   
		&3t_yL^2\eta^2+c_8\gamma^2         
	\end{bmatrix}
\end{align*}
and the last two columns of $\vA$ is given by
	\begin{align*}
		\begin{bmatrix}
			0                  & c_2\gamma & 3c_2\gamma & c_x+c_6\gamma^2 & 3c_7L^2\gamma^2\\
			0         & 0       & c_2\gamma^2   & 0 & c_y+c_7\gamma^2
		\end{bmatrix}^{\T}.
	\end{align*}



\subsubsection{First inequality in \eqref{LIS}}
\begin{equation}\label{ineq:thm1}
	(1-\frac{3}{2}\eta\mu)\epsilon_1+\frac{3\eta L^2}{\mu n}\epsilon_2 \le (1-\frac{1}{2}\eta\mu)\epsilon_1.
\end{equation}
Inequality \eqref{ineq:thm1} holds if 
$\frac{3\kappa^2}{ n}\epsilon_2  \le \epsilon_1$,
where $\kappa:=L/\mu$ is the condition number.
\subsubsection{Second inequality in \eqref{LIS}}
\begin{align}\label{ineq2:thm1}
	\frac{1+\tilde{\rho}^2}{2}\epsilon_2 
	+c_1\frac{\eta^2}{\gamma}L^2\epsilon_3
	+c_2\gamma\epsilon_4
	\leq(1-\frac{1}{2}\eta\mu)\epsilon_2.
\end{align}
Recalling $\frac{1-\tilde{\rho}^2}{2}=\frac{\gamma s(2-\gamma s)}{2}>\frac{\gamma s}{2}$ for $\gamma<\frac{1}{s}$, relation \eqref{ineq2:thm1} holds if 
\begin{align}\label{ineq2:thm2}
	c_1\frac{\eta^2}{\gamma}L^2\epsilon_3
	+c_2\gamma\epsilon_4
	\leq(\frac{s\gamma}{2}-\frac{1}{2}\eta\mu)\epsilon_2.
\end{align}
Dividing $\gamma$ on both sides of \eqref{ineq2:thm2}, we get
\begin{align}\label{ineq2:thm3}
	\frac{\mu\epsilon_2 }{2}\frac{\eta}{\gamma} +c_1L^2\epsilon_3\frac{\eta^2}{\gamma^2}
	+c_2\epsilon_4
	\leq\frac{ s\epsilon_2}{2}.
\end{align}
It is sufficient that
$
	2c_2\epsilon_4\leq\frac{ s\epsilon_2}{2},~
	4\frac{\mu\epsilon_2 }{2}\frac{\eta}{\gamma}\leq\frac{ s\epsilon_2}{2},~
	c_1L^2\epsilon_3\frac{\eta^2}{\gamma^2}\leq\frac{\mu\epsilon_2 }{2}\frac{\eta}{\gamma}.
$
Therefore, if we have  $\eta\leq \min\{\frac{s}{4\mu}\gamma, \frac{\mu\epsilon_2}{2c_1L^2\epsilon_3}\gamma \}$ and $\epsilon_2\geq 2c_1c_2\epsilon_4$, 
then \eqref{ineq2:thm1} can be demonstrated.

\subsubsection{Third inequality in \eqref{LIS}}
\begin{multline}\label{ineq3:thm1}
 c_3nL^2\frac{\eta^2}{\gamma}\epsilon_1
+(c_3L^2\frac{\eta^2}{\gamma}+c_4L^2\gamma)\epsilon_2
+(\frac{1+\tilde{\rho}^2}{2}+\frac{1}{2}c_3\frac{\eta^2}{\gamma})L^2\epsilon_3\\
 +3c_2L^2\gamma\epsilon_4
+c_2\gamma L^2\epsilon_5
\leq(1-\frac{1}{2}\eta\mu)L^2\epsilon_3.
\end{multline}
Dividing $L^2$ on the both side of \eqref{ineq3:thm1}, we have
\begin{multline}\label{ineq3:thm01}
  c_3n\frac{\eta^2}{\gamma}\epsilon_1
	+(c_3\frac{\eta^2}{\gamma}+c_4\gamma)\epsilon_2
	+(\frac{1+\tilde{\rho}^2}{2}+\frac{1}{2}c_3\frac{\eta^2}{\gamma})\epsilon_3\\
+3c_2\gamma\epsilon_4
	+c_2\gamma \epsilon_5
	\leq(1-\frac{1}{2}\eta\mu)\epsilon_3.
\end{multline}
Based on the similar proof in the second inequality, the third inequality holds if 
\begin{align}\label{ineq3:thm2}
m_1\frac{\eta^2}{\gamma^2}+ m_2+ \frac{\mu\epsilon_3}{2}\frac{\eta}{\gamma}\leq \frac{s\epsilon_3}{2},
\end{align}
where 
$
	m_1=:c_3(n\epsilon_1+\epsilon_2+\frac{1}{2}\epsilon_3),~
	m_2=:c_4\epsilon_2+c_2(3\epsilon_4+\epsilon_5).
$ 
It is sufficient that
$
	2m_2\leq\frac{s\epsilon_3}{2},~
	2\mu\epsilon_3\frac{\eta}{\gamma}\leq\frac{s\epsilon_3}{2},~
	m_1\frac{\eta^2}{\gamma^2}\leq\frac{\mu\epsilon_3}{2}\frac{\eta}{\gamma}.	
$
Thus, if there holds $\epsilon_3\geq\frac{4m_2}{s}$ and $\eta\leq \min\{ \frac{s}{4\mu}\gamma, \frac{\mu\epsilon_3}{2m_1}\gamma\}$, we can demonstrate \eqref{ineq3:thm1}.

\subsubsection{Fourth inequality in \eqref{LIS}}
\begin{multline} \label{ineq4:thm1}
2nt_xL^2\eta^2\epsilon_1
	+(2t_xL^2\eta^2
	+c_5\gamma^2)\epsilon_2\\
	+t_x\eta^2 L^2\epsilon_3
	+(c_x+c_6\gamma^2)\epsilon_4
\leq(1-\frac{1}{2}\eta\mu)\epsilon_4.
\end{multline}
It is equivalent that
\begin{multline} \label{ineq4:thm2}
	t_x(2n\epsilon_1+2\epsilon_2+\epsilon_3)L^2\eta^2\\
	+(c_5\epsilon_2+c_6\epsilon_4)\gamma^2
	+\frac{1}{2}\eta\mu\epsilon_4
	\leq(1-c_x)\epsilon_4.
\end{multline}
Inequality \eqref{ineq4:thm2} holds if $\eta\leq\frac{\gamma}{L}$ and $\gamma\leq\min\{1,\frac{1-c_x}{m_3}\epsilon_4\}$,
where 
$
	m_3=:t_x(2n\epsilon_1+2\epsilon_2+\epsilon_3)+c_5\epsilon_2+c_6\epsilon_4
	+\frac{\epsilon_4}{2\kappa}.
$

\subsubsection{Fifth inequality in \eqref{LIS}}
\begin{multline}
	\label{ineq5:thm1}
	6nt_yL^4\eta^2\epsilon_1
	+(6t_yL^4\eta^2+3c_8L^2\gamma^2)\epsilon_2
	+(c_8\gamma^2+3t_yL^2\eta^2)L^2\epsilon_3\\
	+3c_7L^2\gamma^2\epsilon_4
	+(c_y+c_7\gamma^2)L^2\epsilon_5
\leq(1-\frac{1}{2}\eta\mu)L^2\epsilon_5.
\end{multline}
Dividing $L^2$ on both sides of \eqref{ineq5:thm1}, it is equivalent that
\begin{multline}\label{ineq5:thm2}
 	t_y(6n\epsilon_1+6\epsilon_2+3\epsilon_3)L^2\eta^2
	+(3c_8\epsilon_2+c_8\epsilon_3\\
	\qquad+3c_7\epsilon_4+c_7\epsilon_5)\gamma^2
	+\frac{\mu\epsilon_5}{2}\eta 
	\leq(1-t_y)\epsilon_5.
\end{multline}
Inequality \eqref{ineq5:thm2}  holds if $\eta\leq\frac{\gamma}{L}$ and $\gamma\leq\min\{1,\frac{1-c_y}{m_4}\epsilon_5\}$,
where 
$
m_4=:3t_y(2n\epsilon_1+2\epsilon_2+\epsilon_3)
	+3c_8\epsilon_2+c_8\epsilon_3+3c_7\epsilon_4
+c_7\epsilon_5
	+\frac{\epsilon_5}{2\kappa}.
$

In short,  if the positive constants $\epsilon_1$-$\epsilon_5$, consensus step-size $\gamma$ and step-size $\eta$ satisfy the following conditions,	
\begin{align}
	n\epsilon_1 \geq  &3\kappa^2\epsilon_2,
	\epsilon_2\geq 2c_1c_2\epsilon_4,
	\epsilon_3\geq \frac{4m_2}{s},
	\epsilon_4>0,
	\epsilon_5>0,
	\label{condition_epsilon}\\
	\gamma\leq &\min\Big\{1,
	      \frac{1-c_x}{m_3}\epsilon_4,
	      \frac{1-c_y}{m_4}\epsilon_5\Big\},\label{gammaCGT}\\
	\eta\le &\min\Big\{
	     \frac{\mu\epsilon_2}{2c_1L^2\epsilon_3}\gamma,
	     \frac{\mu\epsilon_3}{2m_1}\gamma,
	     \frac{s}{4\mu}\gamma,
	     \frac{\gamma}{L}\Big\},\label{etaCGT}
\end{align}
we can demonstrate the linear system of inequalities in (\ref{LIS}).  
For \eqref{condition_epsilon}, it is easy to verify that there exist solutions to $\epsilon_1-\epsilon_5$. 
Recalling $c_1,c_3$ and $m_1$, we obtain
$\frac{\mu\epsilon_2}{2c_1L^2\epsilon_3}\gamma=\frac{s\epsilon_2}{4\kappa\epsilon_3}\frac{\gamma}{L}$
and $\frac{\mu\epsilon_3}{2m_1}\gamma=\frac{s\epsilon_3}{12\kappa(2n\epsilon_1+2\epsilon_2+\epsilon_3)}\frac{\gamma}{L}$.
Note that $\frac{s\epsilon_3}{12\kappa(2n\epsilon_1+2\epsilon_2+\epsilon_3)}\frac{\gamma}{L}<\frac{s}{4\mu}\gamma$ and $\frac{s\epsilon_3}{12\kappa(2n\epsilon_1+2\epsilon_2+\epsilon_3)}\frac{\gamma}{L}<\frac{\gamma}{L}$,
we have 
\begin{align}	\eta\le &\min\Big\{
	\frac{s\epsilon_2}{4\kappa\epsilon_3}\frac{\gamma}{L},
	\frac{s\epsilon_3}{12\kappa(2n\epsilon_1+2\epsilon_2+\epsilon_3)}\frac{\gamma}{L}
	\Big\}.\label{etaCGT01}
\end{align}

%% file: section/supp/PfEFCGT.tex
\subsection{Proof of Lemma \ref{Lem:EFCGT}}\label{Pf:LemEFCGT}
Before deriving the linear system of inequality, we bound $\EE\left[\|\vX^k-\widehat{\vX}^k\|^2\middle |\cF^k\right]$ and $\EE\left[\|\vY^k-\widehat{\vY}^k\|^2\middle |\cF^k\right]$, respectively. 
From Lines 4 and 5 in Algorithm \ref{Alg:EFCGT}, we know
\begin{align}
	\nonumber  & \EE\left[\norm{\vX^k-\widehat{\vX}^k}^2\middle |\cF^k\right]= \EE\left[\norm{\vE^{k+1}_x-\vE^k_x}^2\middle |\cF^k\right]\\
	\nonumber \leq &  \tau'_1\EE\left[\norm{\vE^{k+1}_x}^2\middle |\cF^k\right]+\frac{\tau'_1}{\tau'_1-1}\norm{\vE^k_x}^2\\ 
	\nonumber \leq &  \tau'_1(1-\delta)\norm{\vE^k_x+\vX^k-\vH^k_x}^2+\frac{\tau'_1}{\tau'_1-1}\norm{\vE^k_x}^2\\ 
	\nonumber \leq &  \tau'_1\left[(1-\frac{\delta}{2})\norm{\vX^k-\vH^k_x}^2+\frac{2(1-\delta)}{\delta}\norm{\vE^k_x}^2\right]\\
	&+\frac{\tau'_1}{\tau'_1-1}\norm{\vE^k_x}^2 .
\end{align}
Taking $\tau'_1=\frac{1}{1-\delta/2}$, we obtain
\begin{align}
	\nonumber  & \EE\left[\norm{\vX^k-\widehat{\vX}^k}^2\middle |\cF^k\right]\\
	\nonumber \leq &  \norm{\vX^k-\vH^k_x}^2+\frac{2}{2-\delta}\cdot\frac{2(1-\delta)}{\delta}\norm{\vE^k_x}^2+\frac{2}{\delta}\norm{\vE^k_x}^2 \\
	\label{XkminusXkhat} \leq & \norm{\vX^k-\vH^k_x}^2+\frac{6}{\delta}\norm{\vE^k_x}^2 .
\end{align}
Similarly, we have
\begin{align}
	\label{YkminusYkhat} \EE\left[\norm{\vY^k-\widehat{\vY}^k}^2\middle |\cF^k\right]
	\leq   \norm{\vY^k-\vH^k_y}^2+\frac{6}{\delta}\norm{\vE^k_y}^2 .
\end{align}

\subsubsection{Optimality error}
Based on the same computation in Appendix \ref{Pf:LemCGT}, we have
$
	\norm{\oX^{k+1} - \vx^*}^2
	\leq(1-\frac{3}{2}\eta\mu) \norm{\oX^k-\vx^*}^2 + \frac{3\eta L^2}{\mu n}\norm{\vX^k-\vone \oX^k}^2.
$

\subsubsection{Consensus error}
Based on the same derivation as \eqref{ineq:AvgConsensus}, we have
\begin{align*}
	&\EE \left[\norm{\vX^{k+1} - \vone \oX^{k+1}}^2\middle|\cF^k\right] \\
	\leq &\frac{1+\tilde{\rho}^2}{2}\norm{\vX^k-\vone\oX^k}^2 +\frac{2\eta^2}{s \gamma}\norm{\vY^k-\vone \oY^k}^2\\
	&+\frac{2\gamma}{s}\norm{\vI-\vW}^2 \EE \left[\norm{\vX^k-\widehat{\vX}^k}^2 \middle| \cF^k\right].
\end{align*}
Thus, recalling \eqref{XkminusXkhat}, we get
\begin{align*}
	&\EE \left[\norm{\vX^{k+1} - \vone \oX^{k+1}}^2\middle|\cF^k\right] \\
	\leq &\frac{1+\tilde{\rho}^2}{2}\norm{\vX^k-\vone\oX^k}^2 +\frac{2\eta^2}{s \gamma}\norm{\vY^k-\vone \oY^k}^2\\
	&+\frac{2\gamma}{s}\norm{\vI-\vW}^2 \norm{\vX^k-\vH_x^k}^2 + \frac{12\gamma}{s\delta}\norm{\vI-\vW}^2\norm{\vE_x^k}^2.
\end{align*}

\subsubsection{Gradient tracker error}

Following from the same computation as \eqref{YkplusCGT}, we get
\begin{align}
	\nonumber &\EE \left[\norm{\vY^{k+1} - \vone \oY^{k+1}}^2\middle|\cF^k\right] \\
	\nonumber= &\frac{1+\tilde{\rho}^2}{2}\norm{\vY^k-\vone\oY^k}^2 +\frac{2}{s\gamma}\EE\bigg[\norm{\nabla^{k+1}-\nabla^k}^2\bigg]\\
	\label{Ykplus}&+\frac{2\gamma}{s}\norm{\vI-\vW}^2\EE\left[\norm{\vY^k-\widehat{\vY}^k}^2 \middle| \cF^k\right]
\end{align}
As before, we obtain
\begin{align*}
	\EE\left[\norm{\nabla^{k+1}-\nabla^k}^2\middle| \cF^k\right]\leq L^2\EE\left[\norm{\vX^{k+1}-\vX^k}^2\middle| \cF^k\right].
\end{align*}
Based on the same derivation as \eqref{ineq:XkplusminusXk}--\eqref{ineq:Ykbar}, we have 
\begin{align}
	\nonumber &\EE\left[\norm{\vX^{k+1}-\vX^k}^2\middle| \cF^k\right]\\
	\nonumber\leq& 3\gamma^2\norm{\vI-\vW}^2\EE\left[\norm{\vX^k-\widehat{\vX}^k}^2\middle| \cF^k\right]\\
	\nonumber&+3\gamma^2\norm{\vI-\vW}^2\norm{\vX^k-\vone\oX^k}^2+3\eta^2\norm{\vY^k-\vone\oY^k}^2\\
	\label{EFineq:XkplusminusXk1}&+6L^2\eta^2\norm{\vX^k-\vone\oX^k}^2+6nL^2\eta^2\norm{\oX^k-\vx^*}^2\\
	\nonumber \leq&  3\gamma^2\norm{\vI-\vW}^2\norm{\vX^k-\vH^k_x}^2+\frac{18\gamma^2}{\delta}\norm{\vI-\vW}^2\norm{\vE^k_x}^2\\
	\nonumber & +6nL^2\eta^2\norm{\oX^k-\vx^*}^2 +3\eta^2\norm{\vY^k-\vone\oY^k}^2\\
	&+(6L^2\eta^2+3\gamma^2\norm{\vI-\vW}^2)\norm{\vX^k-\vone\oX^k}^2
	\label{EFineq:XkplusminusXk2} ,
\end{align}
where the last inequality is from \eqref{XkminusXkhat}. Therefore, we get
\begin{multline}\label{GradplusminusGradEFCGT}
	\EE\left[\norm{\nabla^{k+1}-\nabla^k}^2\middle| \cF^k\right]
	\leq  6nL^4\eta^2\norm{\oX^k-\vx^*}^2\\
	+(6L^4\eta^2+3L^2\gamma^2\norm{\vI-\vW}^2)\norm{\vX^k-\vone\oX^k}^2\\
	~~+3L^2\eta^2\norm{\vY^k-\vone\oY^k}^2+3L^2\gamma^2\norm{\vI-\vW}^2\norm{\vX^k-\vH^k_x}^2\\
	+\frac{18L^2\gamma^2}{\delta}\norm{\vI-\vW}^2\norm{\vE^k_x}^2.
\end{multline}
Plugging \eqref{GradplusminusGradEFCGT} and \eqref{YkminusYkhat} into \eqref{Ykplus}, we conclude that 
\begin{align}
	\nonumber &\EE \left[\norm{\vY^{k+1} - \vone \oY^{k+1}}^2\middle|\cF^k\right] \\
	\nonumber \leq &  \frac{12nL^4\eta^2}{s\gamma}\norm{\oX^k-\vx^*}^2\\
	\nonumber&+(\frac{12L^4\eta^2}{s\gamma}+\frac{6L^2\gamma}{s}\norm{\vI-\vW}^2)\norm{\vX^k-\vone\oX^k}^2\\
	\nonumber&+(\frac{1+\tilde{\rho}^2}{2}+\frac{6L^2\eta^2}{s\gamma})\norm{\vY^k-\vone\oY^k}^2\\
    \nonumber &+\frac{6L^2\gamma}{s}\norm{\vI-\vW}^2\norm{\vX^k-\vH^k_x}^2\\
	\nonumber&+\frac{2\gamma}{s}\norm{\vI-\vW}^2\norm{\vY^k-\vH_y^k}^2
	+\frac{36L^2\gamma}{s\delta}\norm{\vI-\vW}^2\norm{\vE^k_x}^2\\
	\label{Ykplus1EFCGT}
	&+\frac{12\gamma}{s\delta}\norm{\vI-\vW}^2\norm{\vE^k_y}^2.
\end{align}

\subsubsection{Difference compression error of decision variables}
Similarly to the derivation of \eqref{ineq:EHkxError}, we have 
\begin{align}
	\nonumber & \EE \left[\norm{\vX^{k+1} - \vH_x^{k+1}}^2\middle|\cF^k\right]\\
	\leq& d_x\norm{\vX^k-\vH_x^k}^2
	\label{EFineq:EHkxError}	+\frac{\tau'_x}{\tau'_x-1}\EE \left[\norm{\vX^{k+1} -\vX^k}^2\middle|\cF^k\right],
\end{align}
where $d_x=\tau'_x[\alpha_x(1-\delta)+(1-\alpha_x)]<1$.
Plugging \eqref{EFineq:XkplusminusXk2} into \eqref{EFineq:EHkxError}, we obtain
\begin{align}
	\nonumber & \EE \left[\norm{\vX^{k+1} - \vH_x^{k+1}}^2\middle|\cF^k\right]\\
	\nonumber\leq& \frac{\tau'_x}{\tau'_x-1}\Big[6nL^2\eta^2\norm{\oX^k-\vx^*}^2\\
	\nonumber&+(6L^2\eta^2+3\gamma^2\norm{\vI-\vW}^2)\norm{\vX^k-\vone\oX^k}^2\\
	\nonumber&+3\eta^2\norm{\vY^k-\vone\oY^k}^2+\frac{18\gamma^2}{\delta}\norm{\vI-\vW}^2\norm{\vE^k_x}^2\Big]\\
&+\Big[d_x+\frac{3\tau'_x}{\tau'_x-1}\gamma^2\norm{\vI-\vW}^2\Big]\norm{\vX^k-\vH^k_x}^2.
\end{align}

\subsubsection{Difference compression error of gradient trackers}
Similarly to the derivation of \eqref{ineq:EHkxError}, we have 
\begin{align}
	\nonumber & \EE \left[\norm{\vY^{k+1} - \vH_y^{k+1}}^2\middle|\cF^k\right]\\
	\label{EFineq:EHkyError}\leq& d_y\norm{\vY^k-\vH_y^k}^2
+\frac{\tau'_y}{\tau'_y-1}\EE \left[\norm{\vY^{k+1} -\vY^k}^2\middle|\cF^k\right],
\end{align}
where $d_y=\tau'_y\left[\alpha_y(1-\delta)+(1-\alpha_y)\right]<1$. 
Taking the same bound technique as \eqref{ineq:EHkyError1}, we get
\begin{multline}
 \EE \left[\norm{\vY^{k+1} -\vY^k}^2\middle|\cF^k\right]
	\leq3\gamma^2\norm{\vI-\vW}^2\norm{\vY^k-\vone\oY^k}^2\\
+3\gamma^2\norm{\vI-\vW}^2\EE\left[\norm{\vY^k-\widehat{\vY}^k}^2\middle|\cF^k\right]\\
	\label{EFineq:EHkyError1}+3\EE\left[\norm{\nabla^{k+1}-\nabla^k}^2\middle|\cF^k\right].
\end{multline}
Plugging \eqref{YkminusYkhat} and \eqref{GradplusminusGradEFCGT} into \eqref{EFineq:EHkyError1}, we conclude that 
\begin{align}
	\nonumber  &\EE \left[\norm{\vY^{k+1} -\vY^k}^2\middle|\cF^k\right]\\
	\nonumber \leq&
	3\gamma^2\norm{\vI-\vW}^2\norm{\vY^k-\vH^k_y}^2
	+18nL^4\eta^2\norm{\oX^k-\vx^*}^2\\
	\nonumber&+(18L^4\eta^2+9L^2\gamma^2\norm{\vI-\vW}^2)\norm{\vX^k-\vone\oX^k}^2\\
	\nonumber &+(3\gamma^2\norm{\vI-\vW}^2+9L^2\eta^2)\norm{\vY^k-\vone\oY^k}^2\\
	\nonumber&+9L^2\gamma^2\norm{\vI-\vW}^2\norm{\vX^k-\vH^k_x}^2
	 +\frac{54L^2\gamma^2}{\delta}\norm{\vI-\vW}^2\norm{\vE^k_x}^2\\
	&+\frac{18\gamma^2}{\delta}\norm{\vI-\vW}^2\norm{\vE^k_y}^2.
	\label{EFineq:EHkyError2}
\end{align}
Substituting \eqref{EFineq:EHkyError2} into \eqref{EFineq:EHkyError}, we obtain
\begin{align}
	\nonumber & \EE \left[\norm{\vY^{k+1} - \vH_y^{k+1}}^2\middle|\cF^k\right]\\
	\nonumber  \leq& \Big[d_y
+\frac{3\tau'_y}{\tau'_y-1}\gamma^2\norm{\vI-\vW}^2\Big]\norm{\vY^k-\vH_y^k}^2
	\label{EFineq:EHkyError3}\\
	\nonumber  &+\frac{\tau'_y}{\tau'_y-1}\Big[18nL^4\eta^2\norm{\oX^k-\vx^*}^2\\
	\nonumber&+(18L^4\eta^2+9L^2\gamma^2\norm{\vI-\vW}^2)\norm{\vX^k-\vone\oX^k}^2\\
	\nonumber &+(3\gamma^2\norm{\vI-\vW}^2+9L^2\eta^2)\norm{\vY^k-\vone\oY^k}^2\\
	\nonumber&+9L^2\gamma^2\norm{\vI-\vW}^2\norm{\vX^k-\vH^k_x}^2
	+\frac{54L^2\gamma^2}{\delta}\norm{\vI-\vW}^2\norm{\vE^k_x}^2\\
	&+\frac{18\gamma^2}{\delta}\norm{\vI-\vW}^2\norm{\vE^k_y}^2
	\Big].
\end{align}


\subsubsection{Error feedback  of decision variables}
\begin{multline}
 \EE \left[\norm{\vE^{k+1}_x}^2 \middle|\cF^k\right]\leq (1-\delta)\norm{\vE^k_x+\vX^k-\vH^k_x}^2\\ 
	\label{EFineq: Ekx}\leq (1-\frac{\delta}{2})\norm{\vE^k_x}^2+\frac{2(1-\delta)}{\delta}\norm{\vX^k-\vH^k_x}^2.
\end{multline}
\subsubsection{Error feedback  of gradient trackers}
\begin{multline}
\EE \left[\norm{\vE^{k+1}_y}^2 \middle|\cF^k\right]\leq (1-\delta)\norm{\vE^k_y+\vY^k-\vH^k_y}^2\\ 
	\label{EFineq: Eky}\leq (1-\frac{\delta}{2})\norm{\vE^k_y}^2+\frac{2(1-\delta)}{\delta}\norm{\vY^k-\vH^k_y}^2.
\end{multline}

\subsection{Proof of Theorem \ref{Thm:EFCGT}}\label{Pf:ThmEFCGT}
For notational simplicity,  
denote 
$t'_x=\frac{3\tau'_x}{\tau'_x-1}$,
and $t'_y=\frac{3\tau'_y}{\tau'_y-1}$. 
Then, we know$t'_x>1$ and $t'_y>1$. 
Similarly to the proof of Theorem \ref{Thm:CGT}, we consider the following linear system of inequalities:
\begin{align}\label{EFLIS}
	\vB\boldsymbol{\epsilon}<(1-\frac{1}{2}\eta\mu)\boldsymbol{\epsilon},
\end{align}
where $\boldsymbol{\epsilon}:=[\epsilon_1,\epsilon_2,L^2\epsilon_3,\epsilon_4,L^2\epsilon_5,\epsilon_6,L^2\epsilon_7]^{\T}$, 
the submatrix formed by the first three columns of $\vB$ is given by 
\begin{align*}
		\begin{bmatrix}
			1-\frac{3}{2}\eta\mu&\frac{3\eta L^2}{\mu n}&0\\
			0&\frac{1+\tilde{\rho}^2}{2}& d_1\frac{\eta^2}{\gamma} \\
			6nd_1L^2\frac{\eta^2}{\gamma}&3d_2L^2\gamma+6d_1L^4\frac{\eta^2}{\gamma}&\frac{1+\tilde{\rho}^2}{2}+3d_1L^2\frac{\eta^2}{\gamma}\\
			2nt'_xL^2\eta^2&d_3\gamma^2+2t'_xL^2\eta^2&t'_x\eta^2\\
			6nt'_yL^4\eta^2&3t'_yL^2\gamma^2+6t'_yL^4\eta^2&3t'_yL^2\eta^2+d_4\gamma^2\\
			0&0&0\\
			0&0&0
		\end{bmatrix}
\end{align*}
 and the submatrix formed by the last four columns of $\vB$ is given by
\begin{align*}
	\begin{bmatrix}
		&0&0&0&0\\
		&d_2\gamma&0&\frac{6d_2\gamma}{\delta}&0\\
		&3d_2L^2\gamma&d_2\gamma&\frac{18d_2L^2}{\delta}\gamma&\frac{6d_2}{\delta}\gamma\\
		&d_x+d_3\gamma^2&0&\frac{6d_3}{\delta}\gamma^2&0\\
		&3d_4L^2\gamma^2&d_y+d_4\gamma^2&\frac{18d_4}{\delta}L^2\gamma^2&\frac{6d_4}{\delta}\gamma^2\\
		&\frac{2(1-\delta)}{\delta}&0&1-\frac{\delta}{2}&0\\
		&0&\frac{2(1-\delta)}{\delta}&0&1-\frac{\delta}{2}
	\end{bmatrix}.
\end{align*}


\subsubsection{First inequality in \eqref{EFLIS}}
Similarly,  we have
$
\frac{3\kappa^2}{ n}\epsilon_2  \le \epsilon_1.
$
\subsubsection{Second inequality in \eqref{EFLIS}}
\begin{multline}\label{EFineq2:EFthm1}
	\frac{1+\tilde{\rho}^2}{2}\epsilon_2 
	+d_1\frac{\eta^2}{\gamma}L^2\epsilon_3
	+d_2\gamma\epsilon_4
	+ \frac{6d_2}{\delta}\gamma\epsilon_6
	\leq(1-\frac{1}{2}\eta\mu)\epsilon_2,
\end{multline}
where $d_1=\frac{2}{s}$ and $d_2=\frac{2}{s}\norm{\vI-\vW}^2$.
Recalling $\frac{1-\tilde{\rho}^2}{2}=\frac{\gamma s(2-\gamma s)}{2}>\frac{\gamma s}{2}$ for $\gamma\leq 1$, relation \eqref{EFineq2:EFthm1} holds if 
\begin{align}\label{EFineq2:EFthm2}
	d_1\frac{\eta^2}{\gamma}L^2\epsilon_3
	+d_2\gamma\epsilon_4
	+ \frac{6d_2}{\delta}\gamma\epsilon_6
	\leq(\frac{s\gamma}{2}-\frac{1}{2}\eta\mu)\epsilon_2.
\end{align}
Dividing $\gamma$ on both sides of \eqref{EFineq2:EFthm2}, we get
\begin{align}\label{EFineq2:EFthm3}
	\frac{\mu\epsilon_2 }{2}\frac{\eta}{\gamma} +d_1L^2\epsilon_3\frac{\eta^2}{\gamma^2}
	+d_2\epsilon_4
	+ \frac{6d_2\epsilon_6}{\delta}
	\leq\frac{ s\epsilon_2}{2}.
\end{align}
It is sufficient that
$
	4d_2\epsilon_4\leq\frac{ s\epsilon_2}{2},~
	4\frac{6d_2\epsilon_6}{\delta}\leq\frac{ s\epsilon_2}{2},~
	4\frac{\mu\epsilon_2 }{2}\frac{\eta}{\gamma}\leq\frac{ s\epsilon_2}{2},~
	d_1L^2\epsilon_3\frac{\eta^2}{\gamma^2}\leq\frac{\mu\epsilon_2 }{2}\frac{\eta}{\gamma}.
$
Therefore, if we have $\epsilon_2\geq 4d_1d_2\epsilon_4$, $\epsilon_2\geq \frac{24d_1d_2\epsilon_6}{\delta}$, and $\eta\leq \min\{\frac{s}{4\mu}\gamma, \frac{\mu\epsilon_2}{2d_1L^2\epsilon_3}\gamma \}$,
then \eqref{EFineq2:EFthm1} can be demonstrated.

\subsubsection{Third inequality in \eqref{EFLIS}}
\begin{multline}\label{EFineq3:EFthm1}
    6d_1nL^4\frac{\eta^2}{\gamma}\epsilon_1
	+(6d_1L^4\frac{\eta^2}{\gamma}+3d_2L^2\gamma)\epsilon_2\\
	+(\frac{1+\tilde{\rho}^2}{2}+3d_1L^2\frac{\eta^2}{\gamma})L^2\epsilon_3
	+3d_2L^2\gamma\epsilon_4
	+d_2\gamma L^2\epsilon_5\\
	+\frac{18d_2L^2}{\delta}\gamma\epsilon_6
	+\frac{6d_2}{\delta}\gamma L^2\epsilon_7
	\leq(1-\frac{1}{2}\eta\mu)L^2\epsilon_3.
\end{multline}
Based on the similar proof in the second inequality, the third inequality holds if 
\begin{align}\label{EFineq3:EFthm2}
	m'_1 L^2\frac{\eta^2}{\gamma^2}+ m'_2+ \frac{\mu\epsilon_3}{2}\frac{\eta}{\gamma}\leq \frac{s\epsilon_3}{2},
\end{align}
where 
$
	m'_1=:6d_1n\epsilon_1+6d_1\epsilon_2+3d_1\epsilon_3,~
	m'_2=:3d_2\epsilon_2+3d_2\epsilon_4+d_2\epsilon_5+\frac{18d_2}{\delta}\epsilon_6+\frac{6d_2}{\delta}\epsilon_7.
$
It is sufficient that
$
	2m'_2\leq\frac{s\epsilon_3}{2},~
	2\mu\epsilon_3\frac{\eta}{\gamma}\leq\frac{s\epsilon_3}{2},~
	m'_1 L^2\frac{\eta^2}{\gamma^2}\leq\frac{\mu\epsilon_3}{2}\frac{\eta}{\gamma}.	
$
Thus, if there holds $\epsilon_3\geq\frac{4m'_2}{s}$ and $\eta\leq \min\{ \frac{s}{4\mu}\gamma, \frac{\mu\epsilon_3}{2L^2m'_1}\gamma\}$, we can demonstrate \eqref{EFineq3:EFthm1}.

\subsubsection{Fourth inequality in \eqref{EFLIS}}
\begin{multline} \label{EFineq4:EFthm1}
2nt'_xL^2\eta^2\epsilon_1
	+(2t'_xL^2\eta^2+d_3\gamma^2)\epsilon_2
	+t'_x\eta^2L^2\epsilon_3\\
	+\frac{6d_3}{\delta}\gamma^2\epsilon_6
	+(d_x+d_3\gamma^2)\epsilon_4
	\leq(1-\frac{1}{2}\eta\mu)\epsilon_4,
\end{multline}
where $d_3=t'_x\norm{\vI-\vW}^2$. 
It is equivalent that
\begin{multline} 
(2nt'_x\epsilon_1+2t'_x\epsilon_2+t'_x\epsilon_3)L^2\eta^2\\
	+(d_3\epsilon_2+d_3\epsilon_4+\frac{6d_3\epsilon_6}{\delta})\gamma^2
    +\frac{1}{2}\eta\mu\epsilon_4
	\leq(1-d_x)\epsilon_4. \label{EFineq4:EFthm2}
\end{multline}
Inequality \eqref{EFineq4:EFthm2} holds if $\eta\leq\frac{\gamma}{L}$ and $\gamma\leq\min\{1,\frac{1-d_x}{m'_3}\epsilon_4\}$,
where 
$
	m'_3=:2nt'_x\epsilon_1+2t'_x\epsilon_2+t'_x\epsilon_3+d_3\epsilon_2+d_3\epsilon_4+\frac{6d_3\epsilon_6}{\delta}+\frac{\epsilon_4}{2\kappa}.
$

\subsubsection{Fifth inequality in \eqref{EFLIS}}
\begin{multline}\label{EFineq5:EFthm1}
	 6nt'_yL^4\eta^2\epsilon_1
	+(6t'_yL^4\eta^2+3t'_yL^2\gamma^2)\epsilon_2
	+(d_4\gamma^2+3t'_yL^2\eta^2)L^2\epsilon_3\\
	+3d_4L^2\gamma^2\epsilon_4
	+\frac{18d_4}{\delta}L^2\gamma^2\epsilon_6
+\frac{6d_4}{\delta}\gamma^2L^2\epsilon_7\\
	+(d_y+d_4\gamma^2)L^2\epsilon_5
	\leq(1-\frac{1}{2}\eta\mu)L^2\epsilon_5,
\end{multline}
where  $d_4=t'_y\norm{\vI-\vW}^2$. It is equivalent that
\begin{align}\label{EFineq5:EFthm2}
	\nonumber&(6nt'_y\epsilon_1+6t'_y\epsilon_2+3t'_y\epsilon_3)L^2\eta^2
	+(3t'_y\epsilon_2+d_4\epsilon_3+3d_4\epsilon_4\\
	&+d_4\epsilon_5+\frac{18d_4}{\delta}\epsilon_6
	+\frac{6d_4}{\delta}\epsilon_7)\gamma^2
	+\frac{\mu\epsilon_5}{2}\eta 
	\leq(1-d_y)\epsilon_5.
\end{align}
Inequality \eqref{EFineq5:EFthm2}  holds if $\eta\leq\frac{\gamma}{L}$ and $\gamma\leq\min\{1,\frac{1-d_y}{m'_4}\epsilon_5\}$,
where 
$
 m'_4=:6nt'_y\epsilon_1+6t'_y\epsilon_2+3t'_y\epsilon_3+3t'_y\epsilon_2+d_4\epsilon_3+3d_4\epsilon_4
	+d_4\epsilon_5+\frac{18d_4}{\delta}\epsilon_6
		+\frac{6d_4}{\delta}\epsilon_7+\frac{\epsilon_5}{2\kappa}.
$

\subsubsection{Error feedback for decision \eqref{EFLIS}}
\begin{align}\label{EFineq6:EFthm1}
	\frac{2(1-\delta)}{\delta}\epsilon_4+(1-\frac{\delta}{2})\epsilon_6\leq(1-\frac{1}{2}\eta\mu)\epsilon_6.
\end{align}
Thus, if there holds $\epsilon_6\geq\frac{8(1-\delta)\epsilon_4}{\delta^2}$ and $\eta\leq \frac{\delta}{2\mu}$, 
then \eqref{EFineq6:EFthm1} is demonstrated.

\subsubsection{Error feedback for gradient tracking \eqref{EFLIS}}
\begin{align}\label{EFineq7:EFthm1}
	\frac{2(1-\delta)}{\delta}L^2\epsilon_5+(1-\frac{\delta}{2})L^2\epsilon_7\leq(1-\frac{1}{2}\eta\mu)L^2\epsilon_7.
\end{align}
Similarly, if there holds $\epsilon_7\geq\frac{8(1-\delta)\epsilon_5}{\delta^2}$ and $\eta\leq \frac{\delta}{2\mu}$, 
then \eqref{EFineq7:EFthm1} is demonstrated.

In short,  if the positive constants $\epsilon_1$-$\epsilon_7$, consensus step-size $\gamma$ and step-size $\eta$ satisfy the following conditions,	
\begin{align}
	\nonumber \epsilon_1 \geq&  \frac{3\kappa^2}{n}\epsilon_2,
	\epsilon_2\geq 4d_1d_2\epsilon_4,\epsilon_2\geq\frac{24d_1d_2\epsilon_6}{\delta},
	\epsilon_3\geq \frac{4m'_2}{s},\\
	\epsilon_4>&0,
	\epsilon_5>0,
	\epsilon_6\geq\frac{8(1-\delta)\epsilon_4}{\delta^2},
	\epsilon_7\geq\frac{8(1-\delta)\epsilon_5}{\delta^2},
	\label{EFcondition_epsilon}\\
	\gamma\leq &\min\Big\{1, 
	\frac{1-d_x}{m'_3}\epsilon_4,
	\frac{1-d_y}{m'_4}\epsilon_5\Big\},\\
	\eta\le &\min\Big\{\frac{\mu\epsilon_3}{2L^2m'_1}\gamma,
	\frac{\mu\epsilon_2}{2d_1L^2\epsilon_3}\gamma,
	\frac{s}{4\mu}\gamma,
	\frac{\gamma}{L},
	\frac{\delta}{2\mu}\Big\},
\end{align}
we can demonstrate the linear system of inequalities in (\ref{EFLIS}). 
Recalling that $d_1=\frac{2}{s}$ and $m'_1=3d_1(2n\epsilon_1+2\epsilon_2+\epsilon_3)$,
we have $\frac{\mu\epsilon_2}{2d_1L^2\epsilon_3}\gamma=\frac{s\epsilon_2}{4\kappa\epsilon_3}\frac{\gamma}{L}$
and $\frac{\mu\epsilon_3}{2L^2m'_1}\gamma
=\frac{s\epsilon_3}{6\kappa(2n\epsilon_1+2\epsilon_2+\epsilon_3)}\frac{\gamma}{L}$. Since $\mu\leq L$ and $\kappa\geq 1$, we obtain 
$\frac{s\epsilon_3}{6\kappa(2n\epsilon_1+2\epsilon_2+\epsilon_3)}\frac{\gamma}{L}<\frac{\gamma}{L}$ and $\frac{s\epsilon_3}{6\kappa(2n\epsilon_1+2\epsilon_2+\epsilon_3)}\frac{\gamma}{L}<\frac{s}{4\mu}\gamma$.
Therefore, we get
\begin{align}
		\eta\le &\min\Big\{
		\frac{s\epsilon_3}{6\kappa(2n\epsilon_1+2\epsilon_2+\epsilon_3)}\frac{\gamma}{L},
	\frac{s\epsilon_2}{4\kappa\epsilon_3}\frac{\gamma}{L},
	\frac{\delta}{2\mu}\Big\}.
\end{align}

%% file: section/supp/SuppAlg.tex

\subsection{Agents' view}
\label{sec:alg_agents_view}
In the main parts, for the sake of simplicity, we only describe the algorithms in matrix forms. 
Here we further   present C-GT and EF-C-GT algorithms from each agent's perspective.

\begin{figure}[htp]
	\centering
	\begin{minipage}{.99\linewidth}
		\begin{algorithm}[H]
			\caption{C-GT algorithm from agents' view}
			\label{Alg:CGTi}
			\textbf{Input:} stopping time $K$, step-size $\eta$, consensus step-size $\gamma$, scaling parameters $\alpha_x, \alpha_y$, and initial values  $\vx_{i}^0$, $\vh_{x_{i}}^0$, $\vh_{y_{i}}^0$, $\vy_{i}^0 = \nabla f_{i}(\vx_{i}^0)$\\ 
			\noindent\textbf{Output:} $\vx_{i}^K, \vx_{i}^K$
			\begin{algorithmic}[1]
				\For {each agent $i\in\{1,2,\ldots,n\}$}
				\State $\vh^0_{x_{i},w} = \sum_{j\in\mathcal{N}_{i}^{\text{in}}\cup\{i\}} w_{ij}\vh_{x_{j}}^0$
				\State $\vh^0_{y_{i},w} = \sum_{j\in\mathcal{N}_{i}^{\text{in}}\cup\{i\}} w_{ij}\vh_{y_{j}}^0$
				\EndFor				
				\For{$k=0,1,2,\dots, K-1$} in parallel for each agent $i\in\{1,2,\ldots,n\}$
				\State $\vq_{x_{i}}^k=$\textbf{Compress}($\vx_{i}^k-\vh_{x_{i}}^k$)
				\State $\widehat{\vx}_{i}^k=\vh_{x_{i}}^k+\vq_{x_{i}}^k$
				\State $\vh_{x_{i}}^{k+1}=(1-\alpha_{x})\vh_{x_{i}}^k+\alpha_{x}\widehat{\vx}_{i}^k$
				\State $\vq_{y_{i}}^k=$\textbf{Compress}($\vy_{i}^k-\vh_{y_{i}}^k$)
				\State $\widehat{\vy}_{i}^k=\vh_{y_{i}}^k+\vq_{y_{i}}^k$
				\State $\vh_{y_{i}}^{k+1}=(1-\alpha_{y})\vh_{y_{i}}^k+\alpha_{y}\widehat{\vy}_{i}^k$
				\State Send $\vq_{x_{i}}^k,\vq_{y_{i}}^k$ to agent $l\in\mathcal{N}_{i}^{\text{out}}$ and receive $\vq_{x_{j}}^k,\vq_{y_{j}}^k$ from agent $j\in\mathcal{N}_{i}^{\text{in}}$ \hfill $\vartriangleright$ Communication
				\State $\widehat{\vx}_{i,w}^k=\vh_{x_{i},w}^k+\sum_{j\in\mathcal{N}_{i}^{\text{in}}\cup\{i\}}w_{ij}\vq_{x_{j}}^k$
				\State $\widehat{\vy}_{i,w}^k=\vh_{y_{i},w}^k+\sum_{j\in\mathcal{N}_{i}^{\text{in}}\cup\{i\}}w_{ij}\vq_{y_{j}}^k$
				\State $\vh_{x_{i},w}^{k+1} = (1-\alpha_{x})\vh_{x_{i},w}^{k} +\alpha_{x}\widehat{\vx}_{i,w}^k$
				\State $\vh_{y_{i},w}^{k+1} = (1-\alpha_{y})\vh_{y_{i},w}^{k} +\alpha_{y}\widehat{\vy}_{i,w}^k$
				\State $\vx_{i}^{k+1}=\vx_{i}^k-\gamma(\widehat{\vx}_{i}^k-\widehat{\vx}_{i,w}^k)-\eta\vy_{i}^k$   \vfill\hspace*{\fill} $\vartriangleright$ Update decision variable
				\State $\vy_{i}^{k+1}=\vy_{i}^k-\gamma(\widehat{\vy}_{i}^k-\widehat{\vy}_{i,w}^k)+\nabla f_{i}(\vx_{i}^{k+1})-\nabla f_{i}(\vx_{i}^{k})$  \vfill\hspace*{\fill}$\vartriangleright$ Update gradient tracker
				\EndFor
			\end{algorithmic}
		\end{algorithm}
	\end{minipage}
\end{figure}

\begin{figure}[htp]
	\centering
	\begin{minipage}{.99\linewidth}
		\begin{algorithm}[H]
			\caption{EF-C-GT algorithm from   agents' view}
			\label{Alg:EFCGTi}
			\textbf{Input:} stopping time $K$, step-size $\eta$, consensus step-size $\gamma$, scaling parameters $\alpha_x, \alpha_y$, and initial values  $\vx_{i}^0$, $\vh_{x_{i}}^0$, $\vh_{y_{i}}^0$, $\vy_{i}^0 = \nabla f_{i}(\vx_{i}^0)$\\ 
			\noindent\textbf{Output:} $\vx_{i}^K, \vx_{i}^K$
			\begin{algorithmic}[1]
				\For {each agent $i\in\{1,2,\ldots,n\}$}
				\State $\vh^0_{x_{i},w} = \sum_{j\in\mathcal{N}_{i}^{\text{in}}\cup\{i\}} w_{ij}\vh_{x_{j}}^0$
				\State $\vh^0_{y_{i},w} = \sum_{j\in\mathcal{N}_{i}^{\text{in}}\cup\{i\}} w_{ij}\vh_{y_{j}}^0$
				\EndFor				
				\For{$k=0,1,2,\dots, K-1$} in parallel for each agent $i\in\{1,2,\ldots,n\}$
				\State $\widehat{\vq}_{x_{i}}^k=$\textbf{Compress}($\ve_{x_i}^k+\vx_{i}^k-\vh_{x_{i}}^k$)
				\State $\ve_{x_{i}}^{k+1}=\ve_{x_i}^k+\vx_{i}^k-\vh_{x_{i}}^k-\widehat{\vq}_{x_{i}}^k$
				\State $\widehat{\vx}_{i}^k=\vh_{x_{i}}^k+\widehat{\vq}_{x_{i}}^k$
				\State $\vq_{x_{i}}^k=$\textbf{Compress}($\vx_{i}^k-\vh_{x_{i}}^k$)
				\State $\vh_{x_{i}}^{k+1}=\vh_{x_{i}}^k+\alpha_{x}\vq_{x_{i}}^k$
				\State $\widehat{\vq}_{y_{i}}^k=$\textbf{Compress}($\ve_{y_i}^k+\vy_{i}^k-\vh_{y_{i}}^k$)
                \State $\ve_{y_{i}}^{k+1}=\ve_{y_i}^k+\vy_{i}^k-\vh_{y_{i}}^k-\widehat{\vq}_{y_{i}}^k$
                \State $\widehat{\vy}_{i}^k=\vh_{y_{i}}^k+\widehat{\vq}_{y_{i}}^k$
                \State $\vq_{y_{i}}^k=$\textbf{Compress}($\vy_{i}^k-\vh_{y_{i}}^k$)
                 \State $\vh_{y_{i}}^{k+1}=\vh_{y_{i}}^k+\alpha_{y}\vq_{y_{i}}^k$

				\State Send $\vq_{x_{i}}^k,\vq_{y_{i}}^k$, $\widehat{\vq}_{x_{i}}^k$, and $\widehat{\vq}_{y_{i}}^k$ to agent $l\in\mathcal{N}_{i}^{\text{out}}$ and receive $\vq_{x_{j}}^k,\vq_{y_{j}}^k$ , $\widehat{\vq}_{x_{j}}^k$, and $\widehat{\vq}_{y_{j}}^k$ from agent $j\in\mathcal{N}_{i}^{\text{in}}$ \vfill\hspace*{\fill} $\vartriangleright$ Communication
				\State $\widehat{\vx}_{i,w}^k=\vh_{x_{i},w}^k+\sum_{j\in\mathcal{N}_{i}^{\text{in}}\cup\{i\}}w_{ij}\widehat{\vq}_{x_{j}}^k$
				\State $\widehat{\vy}_{i,w}^k=\vh_{y_{i},w}^k+\sum_{j\in\mathcal{N}_{i}^{\text{in}}\cup\{i\}}w_{ij}\widehat{\vq}_{y_{j}}^k$
				\State $\vh_{x_{i},w}^{k+1} = \vh_{x_{i},w}^{k} +\alpha_{x}\sum_{j\in\mathcal{N}_{i}^{\text{in}}\cup\{i\}}w_{ij}\vq_{x_{j}}^k$
				\State $\vh_{y_{i},w}^{k+1}=\vh_{y_{i},w}^{k}
				+\alpha_{y}\sum_{j\in\mathcal{N}_{i}^{\text{in}}\cup\{i\}}w_{ij}\vq_{y_{j}}^k$
				\State $\vx_{i}^{k+1}=\vx_{i}^k-\gamma(\widehat{\vx}_{i}^k-\widehat{\vx}_{i,w}^k)-\eta\vy_{i}^k$   \vfill\hspace*{\fill} $\vartriangleright$ Update decision variable
				\State $\vy_{i}^{k+1}=\vy_{i}^k-\gamma(\widehat{\vy}_{i}^k-\widehat{\vy}_{i,w}^k)+\nabla f_{i}(\vx_{i}^{k+1})-\nabla f_{i}(\vx_{i}^{k})$  \vfill\hspace*{\fill}$\vartriangleright$ Update gradient tracker
				\EndFor
			\end{algorithmic}
		\end{algorithm}
	\end{minipage}
\end{figure}